\documentclass[graybox]{svmult}

\usepackage{type1cm}        
\usepackage{makeidx}         
\usepackage{graphicx}        
\usepackage{multicol}        
\usepackage{multirow}
\usepackage[bottom]{footmisc}

\usepackage{newtxtext}       %
\usepackage{newtxmath}       

\usepackage{latexsym,url,amscd,amsmath}
\usepackage{graphicx}
\usepackage{enumerate}
\usepackage{epstopdf}
\usepackage{algorithm,algorithmic}

\newcommand{\norm}[1]{\left\lVert#1\right\rVert}

\DeclareMathOperator*{\argmin}{argmin}

\newcommand{\shrink}{\mathrm{shrink}}

\newcommand{\TV}{\mathrm{TV}}
\newcommand{\SH}{SH}
\providecommand{\keywords}[1]{\textbf{Keywords: } #1}

\makeindex             

\begin{document}
\title*{Two-stage Geometric Information Guided Image Reconstruction}
\author{Jing Qin \and Weihong Guo}
\institute{
Jing Qin \at Department of Mathematics\\ University of Kentucky\\
Lexington, KY 40506.\\
 \email{jing.qin@uky.edu}
\and
Weihong Guo \at Department of Mathematics, Applied Mathematics and Statistics\\ Case Western Reserve University\\ Cleveland, OH 44106.\\ \email{wxg49@case.edu}
}

\maketitle

\abstract{In compressive sensing, it is challenging to reconstruct image of high quality from very few noisy linear projections. Existing methods mostly work well on piecewise constant images but not so well on piecewise smooth images such as natural images, medical images that contain a lot of details. We propose a two-stage method called GeoCS to recover images with rich geometric information from very limited amount of noisy measurements. The method adopts the shearlet transform that is mathematically proven to be optimal in sparsely representing images containing anisotropic features such as edges, corners, spikes etc. It also uses the weighted total variation (TV) sparsity with spatially variant weights to preserve sharp edges but to reduce the staircase effects of TV. Geometric information extracted from the results of stage I serves as an initial prior for stage II which alternates image reconstruction and geometric information update in a mutually beneficial way. GeoCS has been tested on incomplete spectral Fourier samples. It is applicable to other types of measurements as well. Experimental results on various complicated images show that GeoCS is efficient and generates high-quality images.}


\section{Introduction}
\subsection{Background}\label{sec:background}
Compressive sensing (CS) (cf. the pioneering work \cite{Candes-Tao-06,Donoho-06}) is a new technique in data acquisition realm that reconstructs the signal from fewer measurements than that required by the classical Nyquist-Shannon sampling theorem. This fact makes it very useful in reducing the sensing cost in a variety of applications such as geophysics, astronomy and medical imaging.

The incomplete measurement $b \in \mathbb{C}^m$ of CS is usually linear projections of the underlying image of interest $\bar{u} \in
\mathbb{C}^{n}$ in the form of $b = A\bar{u} + r$ where $r$ is noise. It relies on sparsity/compressibility of $\bar{u}$ itself or under certain transform $\Phi$ to recover it from $b$. When there is no noise ($r = 0$), a straightforward approach to reconstructing $\bar{u}$ is to solve the $L_0$ problem
\begin{equation}\label{eqn:l0}
\min_u \|\Phi u\|_0 \quad\mbox{s.t.}\quad  A u = b.
\end{equation}
However, since the $L_0$ problem is NP-hard, it is common in practice to consider a convex relaxed $L_1$ optimization problem
\begin{equation}\label{eqn:l1}
\min_u \|\Phi u\|_1 \quad\mbox{s.t.}\quad  A u = b,
\end{equation}
which is more computationally efficient. When the sparsifying transform $\Phi$ is orthonormal, the solution of \eqref{eqn:l1} turns out to be the same as that of \eqref{eqn:l0}, and approximate the underlying signal with an overwhelming probability if $A$ satisfies the restricted isometry property (RIP). A commonly used example of orthonormal $\Phi$ is the Haar wavelet transform. More recently, RIP is generalized to D-RIP \cite{Candes-Eldar-Needell-11}, a property that guarantees accurate recovery of images that are nearly sparse in overcomplete/redundant dictionaries. The theoretical results make it more flexible to choose $\Phi$ and reconstruct the signal using $L_1$ optimization \eqref{eqn:l1}.
In the presence of noise, a relaxed model of interest is
\[
\min_u\norm{\Phi u}_1\quad\mbox{s.t.}\quad \norm{Au-b}_2\leq\sigma,
\]
or equivalently
\[
\min_u\lambda \norm{\Phi u}_1 +\frac12\norm{Au-b}_2^2,
\]
where $\lambda>0$ is related to the noise level $\sigma$ in the data. To preserve the discontinuities of $u$ which correspond to the image features, e.g., edges, total variation (TV) is taken into consideration as an additional regularization term \cite{Candes-Romberg-Tao-06,Lustig-Donoho-Pauly-07,Compton-Osher-Bouchard:11}. Then the model with two regularization terms reads as
\begin{equation}
\min_u\beta\norm{u}_{\mathrm{TV}} +\lambda \norm{\Phi u}_1 + \frac12\norm{Au-b}_2^2\label{eqn:TVL1L2}
\end{equation}
where $\Phi$ is the wavelet transform and $\beta$ depends on the characteristics of the underlying image itself.

Wavelet transform and total variation have been used widely in various computer vision and/or imaging science problems. The advantage of wavelet transform is its optimality in approximating signals containing point-wise singularities, but it is widely known that traditional wavelets are not so effective in dealing with singularities in higher dimensions, such as edges in 2D images. TV is optimal in describing piecewise constant images and preserving image edges. As one of the TV based compressive sensing methods, reconstruction from partial Fourier data (RecPF) \cite{Yang-Zhang-Yin-RecPF-10} considers both the TV and the wavelet regularizations. However, it is well known that the TV regularization will cause staircase effects in image restoration \cite{chan2000high} and compressive sensing \cite{guo2014new}. It has been shown that the TV regularization is closely connected to the wavelet one \cite{Cai-Dong-Osher-Shen:2012}. Therefore, the combination of wavelet and TV is not ideal for reconstructing natural images with abundant directional geometric information from few noisy CS measurements.

In this paper, we present a geometric information guided CS (abbreviated as GeoCS) reconstruction method to improve the performance for the situations when the sampling rate is low and/or the noise level is high. The goal is to preserve geometries and fine features with less data required than the state-of-the-art methods. There are three major contributions in our paper.
\begin{itemize}
\item Shearlet transform instead of the widely used wavelet transform is adopted as $\Phi$ to promote the sparsity of signals and thereby reduce the number of measurements required for accurate recovery of them by the fundamental CS theory. Shearlets \cite{Labate-Lim-Kuynoik-Weiss:05,Kutyniok-Guo-Labate:06,Easley-Labate-Lim:08} provide an optimally sparse approximation of piecewise smooth function with $C^2$ singularity curves, e.g., edges, cusps and corners. It combines the power of multiscale methods with the capability of extracting geometry of images.
\item A two-stage method rather than the conventional one-stage method is applied to obtain better recovered images from fewer measurements than the state of the arts. This reconstruction approach adaptively learns the gradient of increasing accuracy to some extent controlling the image geometry. The first stage is to get an initial image reconstruction from CS measurements based on shearlet transform and TV. The second stage starts with geometric information extracted from the result of stage I, and alternates image reconstruction and geometric information update until it converges. Adaptive TV and shearlet transform are used in the second stage. The proposed two stages are significantly different due to the presence of the edge stopping function in the second stage and the corresponding algorithmic differences between the two stages. Massive numerical experiments show that the two-stage approach outperforms the classical one-stage methods.
\item Apply the alternating direction method of multipliers (ADMM) \cite{Glowinski-Marroco:75,Gabay-Mercier:76} or its equivalent split Bregman method \cite{Tom2009} to efficiently solve the optimization problems in both stages. The algorithm achieves fast convergence and produces high-quality images. Convergence of the algorithm at each stage is guaranteed as well.
\end{itemize}

In \cite{Wang-Yin:09}, the binary reweighted $L_1$ regularization is exploited for recovering 1D sparse signals. While it can be applied to recover the sparse wavelet coefficients of an image and hence the image itself, it can hardly take advantages of image edges to improve the recovery. In Edge Guided Compressive Sensing (EdgeCS) \cite{Guo-Yin-10,Guo-Yin-12EdgeCS}, binary edge detection and image reconstruction are performed alternatively in a mutually beneficial way and thus the recovery quality has been improved. In this paper, we consider a more general edge detection whose range is continuous rather than binary. The spatially variant weights associated to the TV ranges between zero and one based on the extracted salient geometric information. In that case, the sharp edges are still able to be preserved while gradual intensity changes in smooth regions can be preserved as well to reduce the staircase effects resulted in TV regularization. In addition, different from EdgeCS as a one-stage method involving TV and wavelets, GeoCS has two stages with TV and shearlets where reliable geometric information from the first stage can be exploited at the second stage to boost performance. More recently, other edge/geometric information guided image reconstruction methods have been proposed, including edge guided CT image reconstruction \cite{cai2014edge,rong2016ct} and edge-guided TV$_p$ regularization for diffuse optical tomography \cite{tong2018edge}.

The paper is organized as follows. We provide a brief review of the shearlet transform in Section \ref{sec:review}, and present the two-stage geometric information guided algorithm in Section \ref{algorithm}. The convergence analysis of the algorithm and practical parameter selection are presented in Section \ref{convergence}. To show the consistent excellence and robustness of the proposed algorithm, plenty of numerical results and comparisons to related work RecPF and EdgeCS are provided in Section \ref{experiments}. Finally, conclusion and remarks are made in Section \ref{conclusion}.

\section{Review of shearlet transform}
\label{sec:review}

The traditional wavelet transform is based on isotropic dilations and thus has limited ability to describe the geometry of multidimensional functions. Directional representation systems such as ridgelets \cite{Candes-DonohoRidgelets:99}, curvelets \cite{Candes-Donoho-Curvelets:00}, contourlets \cite{Do-Vetterli-05}, and shearlets \cite{Kutyniok-Guo-Labate:06} have been designed to provide much more geometric information of multidimensional functions such as images. Curvelets, a tight frame of elongated oscillatory functions at various scales, was first proposed by Cand\`es and Donoho to generalize wavelet. For any $L^2(\mathbb{R}^2)$ function $f$, the $N$ largest term approximation using the curvelet transform has error norm of order $(\log N)^3N^{-2}$. Since the
curvelets are not generated by taking a family of actions on one function as wavelet, it is numerically difficult to implement. In an attempt to provide a better discrete implementation of the curvelets, the contourlet representation is then proposed. It is a discrete time-domain construction, which is designed to achieve essentially the same frequency tiling as the curvelet representation. With the same rate of approximation error decay as curvelets, shearlets have several advantages: efficient implementation, more directional sensitivity and theoretical relation to the multiresolution analysis.

Shearlet transform is an efficient multiscale directional representation of signals, theoretically proven to be optimal up to a log-factor in encoding images with anisotropic features such as edges, corners and other singularities \cite{Labate-Lim-Kuynoik-Weiss:05,Easley-Labate-Lim:08}. Given any function $\psi\in L^2(\mathbb{R}^2)$, the shearlet system is generated by applying the operations of dilation, shear transformation and translation of $\psi$:
\[
\psi_{ast}=|\det M_{as}|^{-\frac12}\psi(M_{as}^{-1}(x-t))
\]
where
\[
M_{as}=\begin{bmatrix}a&-\sqrt{a}s\\0&\sqrt{a}\end{bmatrix}
=\begin{bmatrix}1&-s\\0&1\end{bmatrix}\begin{bmatrix}a&0\\0&\sqrt{a}\end{bmatrix}
:=B_sA_a
\]
with $a\in \mathbb{R}^+, s\in \mathbb{R}, t\in\mathbb{R}^2$, $B_s$ a shear operator and $A_a$ an anisotropic dilation operator. Note in the generation of wavelets, there are only isotropic scaling and translation involved without shearing or anisotropic scaling. The shearlet transform of function $f\in L^2(\mathbb{R}^2)$ is defined as
\[
\SH_\psi(f)(a,s,t)=\langle f,\psi_{ast}\rangle.
\]
The shearlet transform is invertible if the function $\psi$ satisfies the admissibility property
\[
\int_{\mathbb{R}^2} \frac{|\widehat{\psi}(\omega_1,\omega_2)|^2}{|\omega_1|^2} d\omega_1 d\omega_2 <\infty
\]
where $\widehat{\psi}$ is Fourier transform of $\psi$. Given complete shearlet transform coefficients, the original function $f$ can be recovered by
\[
f(x)=\int_{\mathbb{R}^2}\int_{\mathbb{R}}\int_{\mathbb{R}^+}\langle f,\psi_{ast}\rangle\psi_{ast}(x)\frac{da}{a^3}dsdt.
\]
Discrete shearlet transform can be implemented efficiently using the fast Fourier transform. There are three shearlet toolboxes using MATLAB available online: Local Shearlet Toolbox \url{http://www.math.uh.edu/~dlabate/software.html}, Fast Finite Shearlet Transform (FFST)\cite{Hauser2012}, and ShearLab \url{http://www.shearlet.org}. More recently, the shearlet transform has been successfully applied in image processing, e.g., the shearlet-based total variation denoising algorithm \cite{Easley-Labate-Colonna-09}. Because of its higher sparsity of signal representation and ability to capture directional features, the performance of the shearlet transform has also been explored in CS field \cite{Wang-Wang-Hu-Deng:12}.

\section{Proposed model and algorithm}\label{algorithm}
In this section, we present our reconstruction model and analyze how to apply split Bregman to solve the model at each stage. The idea is to use both the shearlet transform and the weighted TV. To enhance the accuracy of weights associated to the TV regularization, we propose a two-stage method. The first stage is to solve a standard $\TV$-$L_1$-$L_2$ model with the shearlet transform to get an initial guess for the underlying image of interest. Note that the extraction of geometry does not work at Stage I since the accuracy is relatively low. In the second stage, we generate the initial spatially variant weights based on the result from stage I, and then alternate image reconstruction and weights update until it converges at this stage. The entire algorithm alternates the two stages until the relative error between two consecutive results is within a tolerance value.

For the shearlet part $\Phi u$, we adopt the FFST algorithm which involves the Fourier transform and the inverse Fourier transform.
Let \[
\Phi u:=\SH(u)= \sum_{i=1}^N\SH_i(u)
\]
where $\SH_i(u)$ is the $i$th subband of shearlet transform of $u$ and $N$ depends on the number of scales in shearlet transform. The $i$th subband of the shearlet transform can be efficiently implemented as componentwise multiplication with a mask matrix denoted by $H_i$ in the frequency space. We have
\[
\SH_i(u)=\mathrm{vec}(\mathcal{F}^{-1}(H_i.*\hat{U}))=\mathrm{vec}(\mathcal{F}^{-1}(H_i)*U):=M_{H_i}u
\]
where $M_{H_i}\in\mathbb{R}^{n^2\times n^2}$, $U$ is the matrix representation of the vectorized image $u$ and $\hat{U}$ is the Fourier transform of $U$.

We demonstrate the idea of the proposed model using partial Fourier sampling, but it can be extended to other linear projection measurements. Let $A=\mathcal{F}_p := PF$ where $P$ is a selection matrix and $F$ is the Fourier transform operator.
For 2D Fourier transform, $F\in\mathbb{R}^{n^2\times n^2}$ is the Kronecker product of two identical $n\times n$ unitary Fourier transform matrices $G$ with
\[
G_{jk}=\frac1{\sqrt{n}}e^{-2\pi\sqrt{-1}(k-1)(j-1)\slash n},\quad j,k=1,\ldots,n.
\]
It can be shown that $F$ satisfies $F^*F=FF^*=I_{n^2}$. By this notation, we get the explicit representation of $M_{H_i}$ as
\[
M_{H_i}=F^*\mathrm{diag}(\mathrm{vec}(H_i))F.
\]
The selection matrix $P\in\mathbb{R}^{k\times n^2}$ is generated simply by deleting the $(n(j-1)+i)$th row of the $n^2\times n^2$ identity matrix if the $(i,j)$th entry of data matrix is not sampled.

\subsection{Stage I: $\mathrm{TV}$-$L_1$-$L_2$ model}
To simplify our discussion, we assume the image to be studied has a square domain. Let $u\in\mathbb{R}^{n^2}$ be the vectorized ground truth image, and $b\in\mathbb{R}^{k}$ ($k\ll n^2$) the given data.
At the first stage, we consider the unconstrained minimization problem with anisotropic discretization of TV as follows:
\begin{equation}\label{prob:stage1orig}
\min_{u\in\mathbb{R}^{n^2}}\beta\sum_{i=1}^2\norm{D_iu}_1+\lambda\sum_{i=1}^N\norm{\SH_i(u)}_1+\frac12\norm{\mathcal{F}_p(u)-b}_2^2
\end{equation}
where $D_1\in\mathbb{R}^{n^2\times n^2}$ ($D_2$) is a horizontal (vertical) first order finite difference operator with periodic boundary conditions.

Due to the non-differentiability of both TV and $L_1$ terms, we introduce auxiliary variables $r_i\in\mathbb{R}^{n^2}$ ($i=1,2$) and $s_i\in\mathbb{R}^{n^2}$ ($i=1,\ldots,N$) such that $r_i=D_iu$ ($i=1,2$) and $s_i=\SH_i(u)$ ($i=1,\ldots,N$) to split the variables. We wish to solve the problem
\begin{equation}\label{prob:stage1orig2}
\begin{aligned}
\min_{u,r_i,s_i} &\beta\sum_{i=1}^2\norm{r_i}_1+\lambda\sum_{i=1}^N\norm{s_i}_1+\frac12\norm{\mathcal{F}_p(u)-b}_2^2
\mbox{s.t.}\quad r_i=D_iu,\,\,\,s_i=\SH_i(u).
\end{aligned}
\end{equation}
After adding the quadratic penalty terms, we get the following unconstrained problem:
\begin{equation}\label{prob:stage1}
\min_{u,r_i,s_i}\beta\sum_{i=1}^2(\norm{r_i}_1+\frac\mu2\norm{r_i-D_iu}_2^2)+\lambda\sum_{i=1}^N(\norm{s_i}_1
+\frac\tau2\norm{s_i-\SH_i(u)}_2^2)+\frac12\norm{\mathcal{F}_p(u)-b}_2^2.
\end{equation}
The above optimization problem is equivalent to \eqref{prob:stage1orig} when $\mu,\tau>0$ go to infinity. Solving \eqref{prob:stage1} using the continuation scheme \cite{Hale2007} is a straightforward method, but it is slow and leads to the ill conditioning of the problem when $\mu,\tau$ are sufficiently large. We hereby apply the split Bregman, which provides fast convergence while the values of $\mu,\tau$ can be fixed. The split Bregman formulation is
\[
\begin{aligned}
\min_{u,r_i,s_i}&\beta\sum_{i=1}^2(\norm{r_i}_1+\frac\mu2\norm{r_i-D_iu-v_i}_2^2)+\frac12\norm{\mathcal{F}_p(u)-b}_2^2\\
&+\lambda\sum_{i=1}^N(\norm{s_i}_1+\frac\tau2\norm{s_i-\SH_i(u)-t_i}_2^2)
\end{aligned}
\]
where $v_i$'s, $t_i$'s are updated by Bregman iterations
\[
\left\{\begin{aligned}
v_i&\leftarrow v_i+\gamma(D_iu-r_i),\quad i=1,2\\
t_i&\leftarrow t_i+\gamma(\SH_i(u)-s_i),\quad i=1,\ldots,N,
\end{aligned}
\right.
\]
with $\gamma>0$ a parameter to be discussed later.

We finally decompose it into three sets of subproblems and apply alternating minimization scheme to get a minimizer iteratively.
\[
\left\{
\begin{aligned}
&\min_{r_i}\norm{r_i}_1+\frac\mu2\norm{r_i-D_iu-v_i}_2^2,\quad i=1,2\\
&\min_{s_i}\norm{s_i}_1+\frac\tau2\norm{s_i-\SH_i(u)-t_i}_2^2,\quad i=1,\ldots,N\\
&\min_{u}\frac{\beta\mu}2\sum_{i=1}^2\norm{r_i-D_iu-v_i}_2^2+\frac12\norm{\mathcal{F}_p(u)-b}_2^2
+\frac{\lambda\tau}2\sum_{i=1}^N\norm{s_i-\SH_i(u)-t_i}_2^2.
\end{aligned}
\right.
\]

The first two subproblems are both in the form of $L_1$-$L_2$ optimization
\[
\min_{x\in\mathbb{R}^m}\delta\norm{x}_1+\frac12\norm{x-v}_2^2,\quad (\delta>0)
\]
whose solution is given by using the shrinkage operator
\[
x=\mathrm{shrink}(v,\delta):=\mathrm{sgn}(v).*\max\{|v|-\delta,0\}
\]
where $\mathrm{sgn}(x):\mathbb{R}^m\rightarrow\mathbb{R}^m$ is componentwise sign function and $.*$ is componentwise multiplication. Then by the similar derivations, the first two subproblems have closed-form solutions using shrinkage.

To solve the last least square subproblem, we consider the corresponding normal equation
\[\begin{aligned}
&\beta\mu\sum_{i=1}^2D_i^T(D_iu-r_i+v_i)+\lambda\tau\sum_{i=1}^NM_{H_i}^*(M_{H_i}u-s_i+t_i)+(PF)^*(PFu-b)=0.
\end{aligned}\]
To circumvent the expensive computation of the inverse matrix, we multiply both sides by $F$ and simply the solution due to the fact that $F$ is unitary. By simplification, the solution is explicitly represented as
\[\begin{aligned}
&\Big(\beta\mu\sum_{i=1}^2F(D_i^TD_i)F^*+\lambda\tau\sum_{i=1}^NF(M_{H_i}^*M_{H_i})F^*+P^*P\Big)Fu\\
&=\beta\mu\sum_{i=1}^2F(r_i-v_i)+\lambda\tau\sum_{i=1}^NF(s_i-t_i)+P^*b.
\end{aligned}\]
Here $P^*P :=\mathrm{diag}(\widetilde{P})$ is a $n^2\times n^2$ diagonal matrix with diagonal value 0, 1 corresponding to nonsampled and sampled entries, respectively. Denoting $P^*b=\tilde{b}$,
the solution can be further written in terms of Fourier transform
\begin{equation}\label{eqn:subpro3}
\begin{aligned}
u&=\mathcal{F}^{-1}\Big((\beta\mu\sum_{i=1}^2\mathcal{F}(r_i-v_i)+\lambda\tau\sum_{i=1}^N\mathcal{F}(s_i-t_i)+\tilde{b})\\
&./(\beta\mu\sum_{i=1}^2\mathrm{diag}(FD_i^TD_iF^*)+\lambda\tau\sum_{i=1}^N\mathrm{diag}(FM_{H_i}^*M_{H_i}F^*)+\widetilde{P})\Big).
\end{aligned}
\end{equation}
Here $./$ means the componentwise division. One more remark about this approach is that since $D_i$'s and $M_{H_i}$'s are circulant matrices which can be diagonalized under the Fourier transform, both $FD_i^TD_iF$ and $FM_{H_i}^*M_{H_i}F^*$ are diagonal matrices. We follow the convention that $0/0 = 0$. The above analysis yields the following algorithm.
\begin{algorithm}
\caption{GeoCS Stage I (solving \eqref{prob:stage1orig2})}
\label{alg1}
\begin{algorithmic}
\item[1.] Initialization: set $u^0,\,r_i^0,v_i^0,s_i^0,t_i^0$ as zero matrices, and choose proper parameters $\beta,\mu,\lambda,\tau,\gamma>0$.
\item[2.] For $k=0,1,2,\ldots$, run the following steps:
\begin{align*}
&r_i^{k+1}=\shrink(D_iu^k+v_i^k,1\slash\mu),\quad i=1,2\\
&s_i^{k+1}=\shrink(\SH_i(u^k)+t_i^k,1\slash\tau),\quad i=1,\ldots,N\\
&u^{k+1}\mbox{ is given by \eqref{eqn:subpro3}}\\
&v_i^{k+1}=v_i^k+\gamma(D_iu^{k+1}-r_i^{k+1}),\quad i=1,2\\
&t_i^{k+1}=t_i^k+\gamma(\SH_i(u^{k+1})-s_i^{k+1}),\quad i=1,2,\ldots,N.\\
&\mbox{If } \|{u^{k+1}-u^k}\|/\|{u^{k+1}}\|\leq \mathrm{tol}, \mbox{ stop the iteration}.
\end{align*}
\end{algorithmic}
\end{algorithm}

\subsection{Stage II: $w\mathrm{TV}$-$L_1$-$L_2$ model}
The stage I model works well in mild CS scenario but not so efficient in challenging scenarios when sampling rate is extremely low and noise is excessive. To handle challenging situations, we start with the result of stage I and then alternatively perform geometric information update and image reconstruction in a beneficial way. Specifically, setting the result from stage I as initial guess, we define adaptive weights based on it, and use weighted TV along with shearlet to reconstruct an image. We then continue alternating weight update and image reconstruction until it converges. For a fixed weight, the model reads as below
\begin{equation}\label{prob:stage2}
\begin{aligned}
\min_{u\in\mathbb{R}^{n^2}}&\beta\sum_{i=1}^2\norm{w_i.*D_iu}_1+\lambda\sum_{i=1}^N\norm{\SH_i(u)}_1+\frac12\norm{\mathcal{F}_p(u)-b}_2^2
\end{aligned}
\end{equation}
where $w_i$'s are the weights based on the extracted geometric information, such as reliable gradients and high frequency subbands of shearlet transform coefficients. And $.*$ is componentwise multiplication. Algorithm \ref{alg2} is designed to address the above problem and the complete stage II algorithm is shown in Algorithm \ref{alg3}.

\begin{algorithm}
\caption{Weighted TV shearlet based image reconstruction algorithm (solving \eqref{prob:stage2})}
\label{alg2}
\begin{algorithmic}
\item[1.] Initialization: set $u^0,\,r_i^0,v_i^0,s_i^0,t_i^0$ as those generated from stage I.
\item[2.] For $k=0,1,2,\ldots$, run the steps:
\begin{align*}
&r_i^{k+1}=\shrink(D_iu^{k}+v_i^{k},w_i\slash\mu),\quad i=1,2\\
&s_i^{k+1}=\shrink(\SH_i(u^{k})+t_i^{k},1\slash\tau),\quad i=1,\ldots,N\\
&u^{k+1}\mbox{ is given by \eqref{eqn:subpro3}}\\
&v_i^{k+1}=v_i^k+\gamma(D_iu^{k+1}-r_i^{k+1}),\quad i=1,2\\
&t_i^{k+1}=t_i^k+\gamma(\SH_i(u^{k+1})-s_i^{k+1}),\quad i=1,2,\ldots,N.\\
&\mbox{If }\|u^{k+1}-u^k\|/\|u^{k+1}\|\leq\mathrm{tol},\mbox{ stop the iteration.}
\end{align*}
\end{algorithmic}
\end{algorithm}

Given the latest iterate $\tilde{u}$ for the reconstructed image, we define TV weight at each pixel as a function of the gradient of $\tilde{u}$ at the same pixel. Suppose $g:[0,\infty)\rightarrow[0,1]$ is a non-increasing function satisfying
\[
g(0)=1,\quad \lim_{s\rightarrow\infty}g(s)=0.
\]
$g$ is called edge stopping function in image segmentation or diffusivity function in PDE. The reason is that $g(|\nabla \tilde{u}|)$ approaches to zero near edges where the gradient gets large while close to one in smooth areas where the gradient becomes small. In fact, besides separating edges from smooth areas, $g$ also identifies the small differences in intensity variations \textit{within} the smooth areas. The pixel in regions of small intensity variations will get larger $g$ value than that in regions of large intensity variations. Weighted TV with this type of weight will preserve the various intensity variation scales in the reconstruction process, and thus increase the robustness of TV and reduce the staircase effects of TV.

\begin{figure}
\begin{center}
\includegraphics[width=0.6\textwidth]{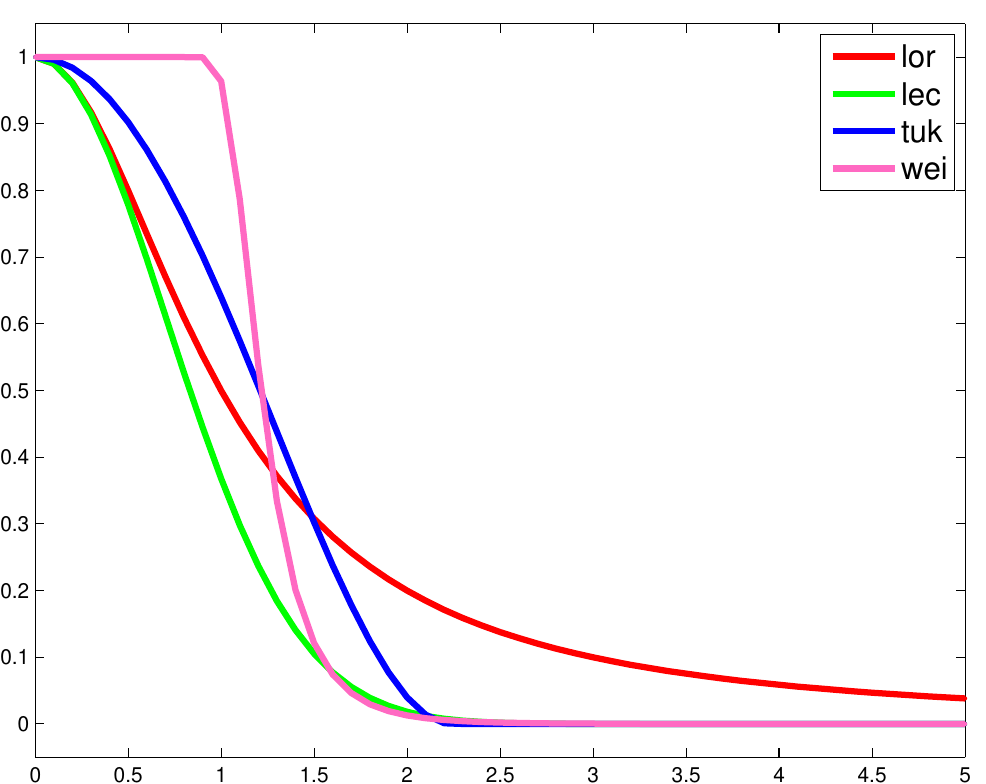}
\end{center}
\caption{Edge stopping function $g$'s when $h=1$}
\label{fig:edgestopfcn}
\end{figure}

There are many choices for $g$. Some commonly used ones are listed below, where $h$ is a parameter controlling the differentiation of smoothness levels.
\begin{enumerate}[(a)]
\item Lorentzian function
\[
g_{_{Lor}}(x)=\frac1{1+\frac{x^2}{h^2}}
\]

\item Le Clerc function
\[
g_{_{Lec}}(x)=\exp\Big(-\frac{x^2}{h^2}\Big)
\]

\item Tukey bi-weight function
\[
g_{_{Tuk}}(x)=\left\{\begin{aligned}
&\Big(1-\frac{x^2}{5h^2}\Big)^2&&|x|<\sqrt{5}h\\
&0&&\mbox{otherwise}
\end{aligned}\right.
\]

\item Weickert function
\[
g_{_{Wei}}(x)=\left\{\begin{aligned}
&1-\exp\Big(-\frac{3.31488h^8}{x^8}\Big)&&x\neq0\\
&1&&\mbox{otherwise}
\end{aligned}\right.
\]
\end{enumerate}
In Fig. \ref{fig:edgestopfcn}, we plot the above four $g$ functions when $h = 1$. From observation, it's clear that they have different decay behaviors. Especially, Weikert edge function decays slowly at the two ends but fast near the middle and Tukey bi-weight function decays slowly all the way long. So for piecewise constant images whose intensity changes sharply from one region to another, Weikert is optimal while Tukey bi-weight function is more appropriate for generic complicated piecewise smooth images with ubiquitous unprecedented intensity variations. We use Tukey bi-weight for all our numerical experiments as we focus on testing piecewise smooth images.

Since we adopt anisotropic TV discretization, our weights are different along $x,y$ directions and are defined respectively as
\begin{equation}\label{eqn:weight}
w_1=g(|D_1\tilde{u}|),\quad w_2=g(|D_2\tilde{u}|).
\end{equation}
Notice that high frequency components of the shearlet transform of $\tilde{u}$ also provide some edge information. So another option to define weights is to gather all the high frequency subbands. But our massive numerical experiments show that $g$ function of gradients is more efficient.

In hope of retrieving more trustworthy geometric prior information, we update the weights from each convergent intermediate result and reapply the Algorithm \ref{alg2}. Then we get the Algorithm \ref{alg3}.
\begin{algorithm}
\caption{GeoCS Stage II}
\label{alg3}
\begin{algorithmic}
\item[1.] Initialization: set $u^0,\,r_i^0,v_i^0,s_i^0,t_i^0$ as those produced from stage I.
\item[2.] For $j=0,1,2,\ldots$, run the steps:
\begin{enumerate}[(1)]
\item Build the weights $w_1^{j+1},\,w_2^{j+1}$ based on $u^j$ by \eqref{eqn:weight}.
\item Set $u^j,r_i^j,v_i^j,s_i^j,t_i^j$ as initial values, apply Algorithm \ref{alg2} to solve \eqref{prob:stage2} and get $u^{j+1},r_i^{j+1},v_i^{j+1},s_i^{j+1},t_i^{j+1}$.
\item If $\|u^{j+1}-u^j\|/\|u^{j+1}\|\leq\mathrm{tol}$, stop the iteration.
\end{enumerate}
\end{algorithmic}
\end{algorithm}

\section{Convergence analysis}\label{convergence}
There are close relationships between Bregman iterative methods and its variants such as linearized Bregman, Bregman operator splitting, and the classical Lagrangian based methods, such as method of multipliers, the alternating direction method of multipliers (ADMM) and alternating minimization algorithm (AMA). The connection between split Bregman algorithm and ADMM, and its illustrative applications in TV-$L_1$ and TV-$L_2$ problems can be found in \cite{Esser2009}. In this section, we aim to bridge the gap between our proposed algorithms using split Bregman method and ADMM by constructing one augmented Lagrangian adapted to our problem. Then the existing convergence theory for ADMM can be used to justify our proposed algorithm utilizing the split Bregman and quadratic penalties.

We first analyze the algorithm in stage I, and the discussions can be analogously extended to the algorithm \ref{alg2} in stage II. Based on the problem \eqref{prob:stage1orig}, we build the augmented Lagrangian as below:
\[\begin{aligned}
&L(u,r_1,r_2,s_1,\ldots,s_N,v_1,v_2,t_1,\ldots,t_N)\\
&=\beta\sum_{i=1}^2(\norm{r_i}_1-v_i^T(r_i-D_iu)+\frac{\mu}2\norm{r_i-D_iu}_2^2)\\
&+\lambda\sum_{i=1}^N(\norm{s_i}_1-t_i^T(s_i-\SH_i(u))+\frac{\tau}2\norm{s_i-\SH_i(u)}_2^2)
+\frac12\norm{\mathcal{F}_p(u)-b}_2^2.
\end{aligned}\]
Since the variables $r_i$'s and $s_i$'s are separable in the Lagrangian $L$, minimizing $L$ over $(r_1,r_2,s_1,\ldots,s_N)$ simultaneously can be replaced by minimizing $L$ over $r_i$'s and $s_i$'s individually. Thus ADMM yields the following iterations
\[
\left\{
\begin{aligned}
r_i^{k+1}&=\argmin_{r_i}\norm{r_i}_1-(v_i^k)^T(r_i^k-D_iu^k)+\frac\mu2\|r_i-D_iu^k\|_2^2,\quad i=1,2\\
s_i^{k+1}&=\argmin_{s_i}\norm{s_i}_1-(t_i^k)^T(s_i-\SH_i(u^k))+\frac\tau2\|s_i-\SH_i(u^k)\|_2^2,\quad i=1,\ldots,N\\
u^{k+1}&=\argmin_{u}L(u,r_1^{k+1},r_2^{k+1},s_1^{k+1},\ldots,s_N^{k+1},v_1^k,v_2^k,t_1^k,\ldots,t_N^k)\\
v_i^{k+1}&=v_i^k+\mu\gamma(D_iu^{k+1}-r_i^{k+1}),\quad i=1,2\\
t_i^{k+1}&=t_i^k+\tau\gamma(\SH_i(u^k)-s_i^{k+1}),\quad i=1,\ldots,N
\end{aligned}
\right.
\]
By absorbing the linear terms involving $v_i$'s and $t_i$'s into the quadratic terms, it can be simplified as
\begin{equation}\label{eqn:ADMM}
\left\{
\begin{aligned}
r_i^{k+1}&=\argmin_{r_i}\norm{r_i}_1+\frac\mu2\Big\|r_i-D_iu^k-\frac{v_i^k}{\mu}\Big\|_2^2,\quad i=1,2\\
s_i^{k+1}&=\argmin_{s_i}\norm{s_i}_1+\frac\tau2\Big\|s_i-\SH_i(u^k)-\frac{t_i^k}{\tau}\Big\|_2^2,\quad i=1,\ldots,N\\
u^{k+1}&=\argmin_{u}\frac{\beta\mu}2\sum_{i=1}^2\|r_i^{k+1}-D_iu\|_2^2
+\frac12\norm{\mathcal{F}_p(u)-b}_2^2+\frac{\lambda\tau}2\sum_{i=1}^N\|s_i^{k+1}-\SH_i(u)\|_2^2\\
v_i^{k+1}&=v_i^k+\mu\gamma(D_iu^{k+1}-r_i^{k+1}),\quad i=1,2\\
t_i^{k+1}&=t_i^k+\tau\gamma(\SH_i(u^k)-s_i^{k+1}),\quad i=1,\ldots,N
\end{aligned}
\right.
\end{equation}
One can see the algorithm derived by ADMM here is equivalent to Algorithm \ref{alg1} by split Bregman method with quadratic penalization. The detailed convergence analysis of the algorithm \eqref{eqn:ADMM} can be found in \cite{Esser2009}. To be complete, we present the convergence theorem without the proof.
\begin{theorem}\label{thm1}
For any $\mu,\tau>0$ and $\gamma\in(0,(\sqrt{5}+1)\slash2)$, the sequences $\{(u^k,r_i^k,s_i^k)\}$ generated by \eqref{eqn:ADMM} from any starting point $(u^0,\lambda^0,\eta^0)$ converges to a solution of problem \eqref{prob:stage1orig2}.
\end{theorem}

Therefore, by choosing an appropriate parameter $\gamma$, the proposed Algorithm \ref{alg1} provides a convergent solution to the problem \eqref{prob:stage1orig}. Likewise for the fixed weights $w_i$'s, by replacing $D_iu^k$ with $w_i.*D_iu^k$ the ADMM yields a similar algorithm equivalent to Algorithm \ref{alg2} and thereby the convergence is guaranteed as well.

\textbf{Parameter selection.} The above theorem only requires ${\gamma\in(0,(\sqrt{5}+1)\slash2)}$, and positive $\mu,\tau$ to guarantee convergence. In the perspective of convergence speed, our experience with a variety of tests shows that $\gamma$ restricted in $(1,\, (\sqrt{5}+1)\slash2)$ consistently yields good results. Refer to \cite{glowinski2008lectures,fortin1983chapter} for using variational techniques to derive the condition on $\gamma$ for the convergence of ADMM. Regarding $\mu$ and $\tau$, as they show up in the shrinkage representation of updates for $r_i$ and $s_i$, an inappropriate selection of them leads to slow convergence. Especially, if $\mu$ and $\tau$ are set too small, the updates for $r_i$ and $s_i$ will dwell in $0$ at the first several iterations. We scale image intensity to $[0,1]$ to make the effect of $\mu$ and $\tau$ on convergence speed moderate. $\beta$ and $\lambda$ depend on the gradient/shearlet transform sparsity of the underlying image and the noise/error level in the measurements. Implementation details and specific parameter selections will be explained in Section \ref{experiments}.

\section{Numerical examples}\label{experiments}
In this section, we illustrate the performance of GeoCS on various images with different sampling rates and noise levels. We also compare GeoCS with two related CS reconstruction approaches: RecPF \cite{Yang-Zhang-Yin-RecPF-10} and EdgeCS \cite{Guo-Yin-10}. All experiments were performed under Windows 7 Professional operating system and MATLAB R2012a running on a Dell desktop with Intel Core i5 CPU at 3.10 GHz and 8 GB of memory.

RecPF iteratively recovers an image from its incomplete Fourier samples by solving
\begin{equation}
\label{RecPF}
\min_u \beta \mathrm{TV}(u) + \lambda \|\Phi u\|_1 + \frac{1}{2}\|\mathcal{F}_p (u) - b\|^2_2
\end{equation}
where $\mathrm{TV}(u)$ can be either isotropic or anisotropic and $\Phi$ is wavelet transform.

EdgeCS alternatively performs image reconstruction and edge detection in a mutually beneficial manner. It detects edges from the intermediate reconstruction and use edge information to guide the next stage of image reconstruction and so on.
GeoCS is different from EdgeCS as analyzed in Section \ref{sec:background}.

Our test images are all piecewise smooth images with a lot of fine details: a human brain MR image, Barbara image with textures and a human knee MR image. The intensity value of each test image is scaled to the range $[0,1]$ before simulating $b$. Partial Fourier CS data are simulated through fast Fourier transform (FFT) on the test images followed by sampling on smooth radial trajectories that are empirically shown to be effective.

All the quantitative comparisons are based on relative error and signal-to-noise ratio (SNR). Relative error is to measure the recovery accuracy and defined as
\[\mathrm{RelErr}=\frac{\norm{u-u_{\mathrm{true}}}_2^2}{\norm{u_{\mathrm{true}}}_2^2}\] where $u$ and $u_{\mathrm{true}}$ are the recovered image and the ground truth, respectively. Considering the independence with the above measure, we adopt the SNR defined in \cite{Coupe08}
\[
\mathrm{SNR}=10\log_{10}\frac{\norm{u.^2+u_{\mathrm{true}}.^2}_2^2}{\norm{u-u_{\mathrm{true}}}_2^2},
\]
where ($.^2$) represents the componentwise squaring. For all the experiments, we fix $\gamma = 1,\mu= \tau = 10^2$ and vary $\beta, \lambda$ slightly based on the noise level. When there is no noise we set $\beta=\lambda=10^{-5}$ for all of the three images while they are set a little larger in the presence of noise. The results are not sensitive to the selection of $\beta$ and $\lambda$. For discrete shearlet transform, we adopt FFST \cite{Hauser2012} with 3 scales and 13 subbands (12 high frequency and one low frequency). The parameter $h$ used in Tukey bi-weight $g$ function in stage II is set among $(0,1]$. The tolerance value in all algorithms is set as $\mathrm{tol}=10^{-5}$. And in each of the following experiment, the total number of iterations used in stage I is less than 1000 and stage II takes less than 100 iterations to achieve a convergent solution.

\subsection{Example 1}
\begin{figure*}
\begin{center}
\begin{tabular}{cccc}
\includegraphics[width=.23\textwidth]{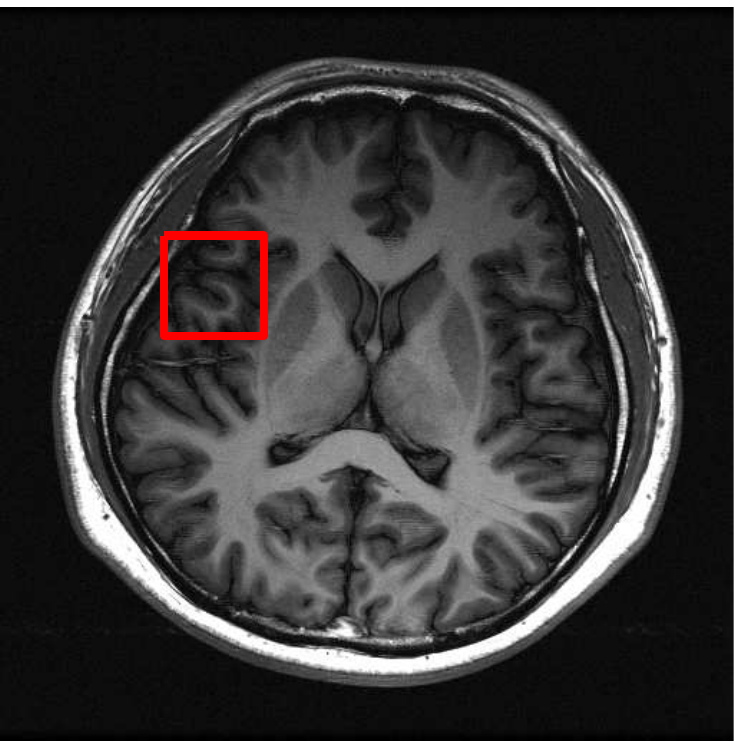}&
\includegraphics[width=.23\textwidth]{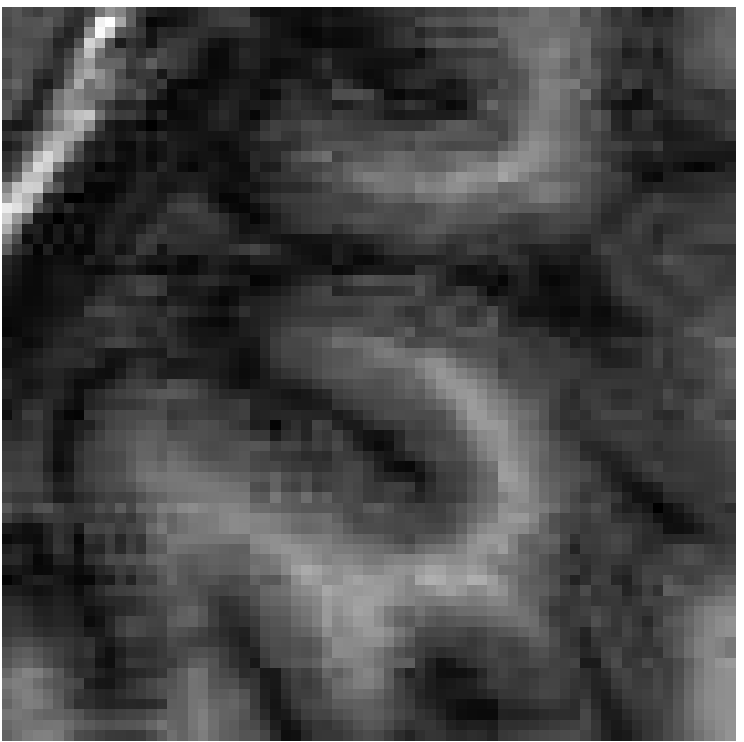}&
\includegraphics[width=.23\textwidth]{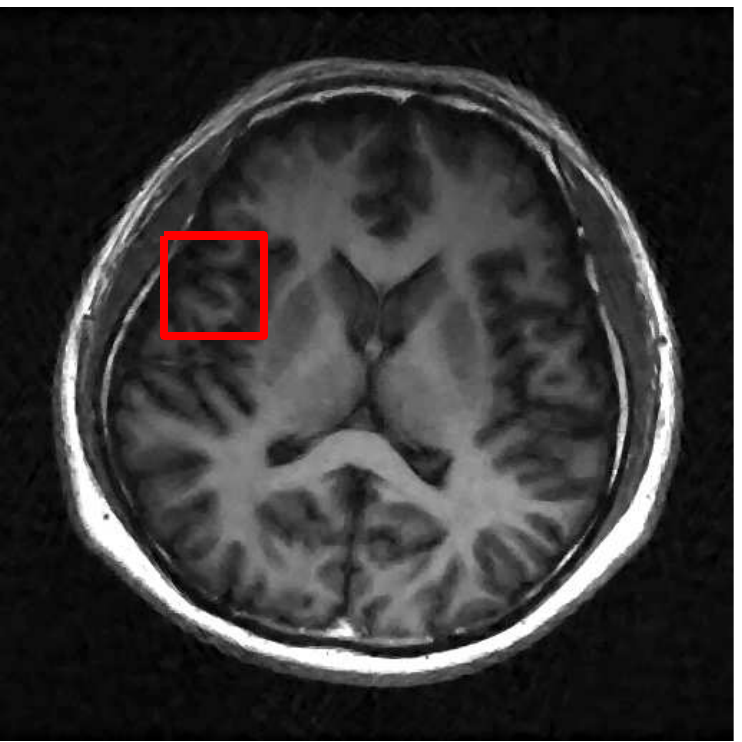}&
\includegraphics[width=.23\textwidth]{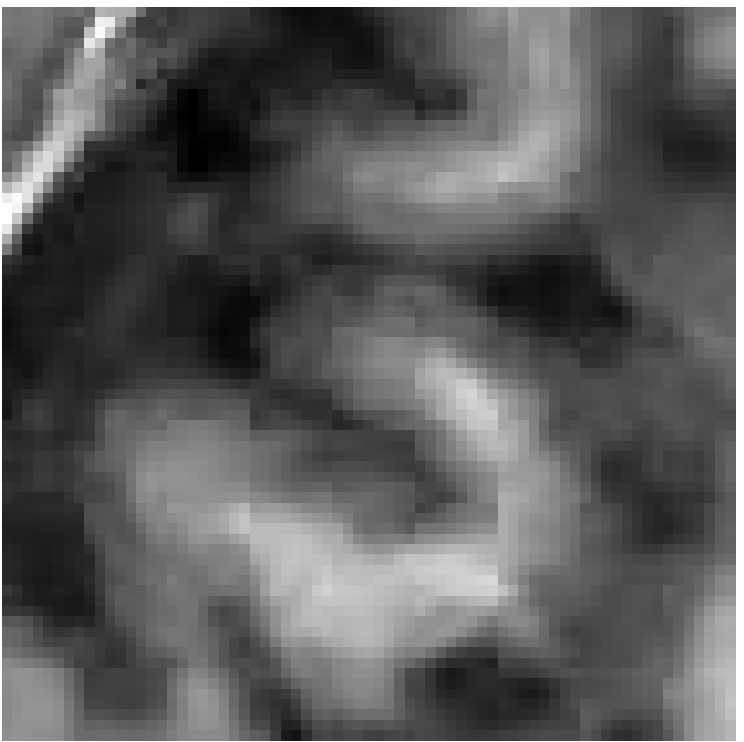}\\
\includegraphics[width=.23\textwidth]{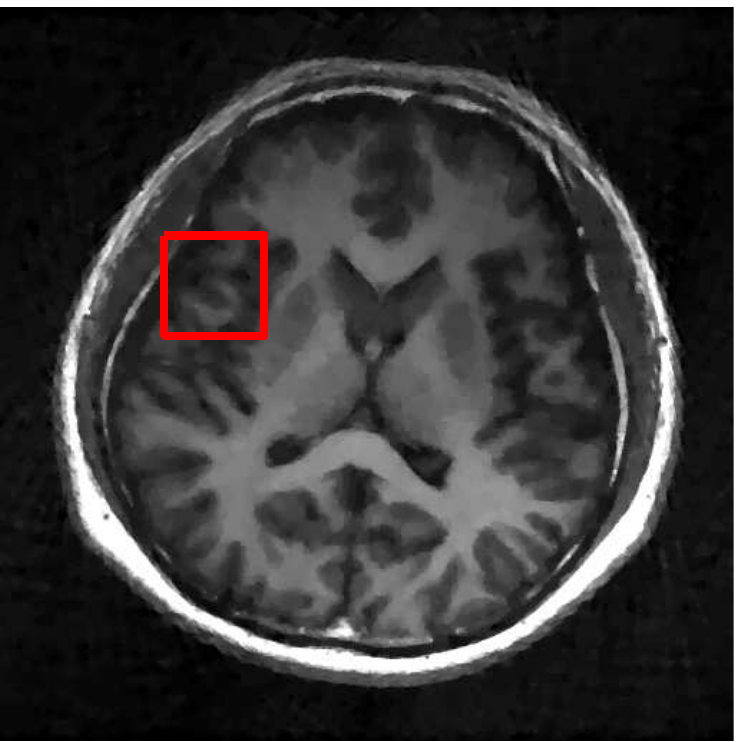}&
\includegraphics[width=.23\textwidth]{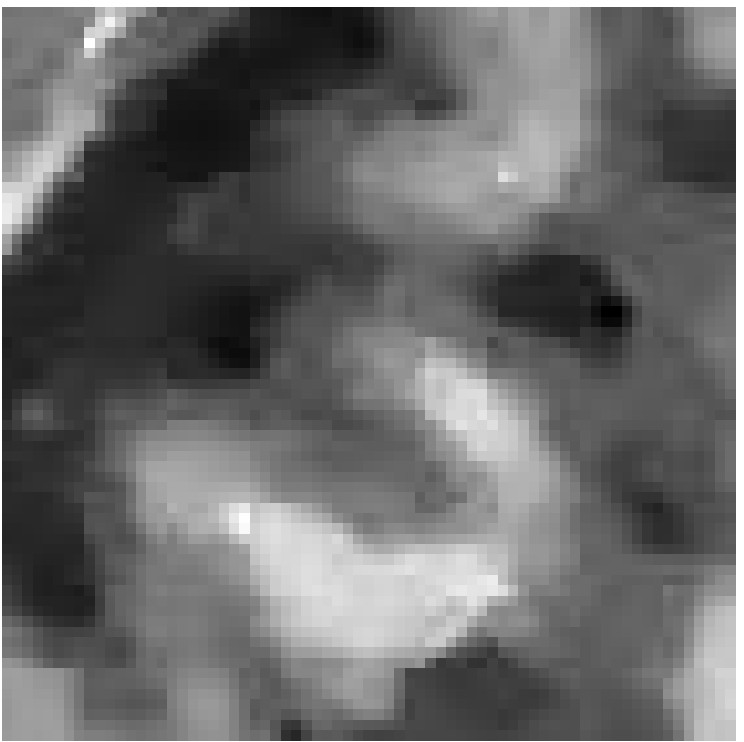}&
\includegraphics[width=.23\textwidth]{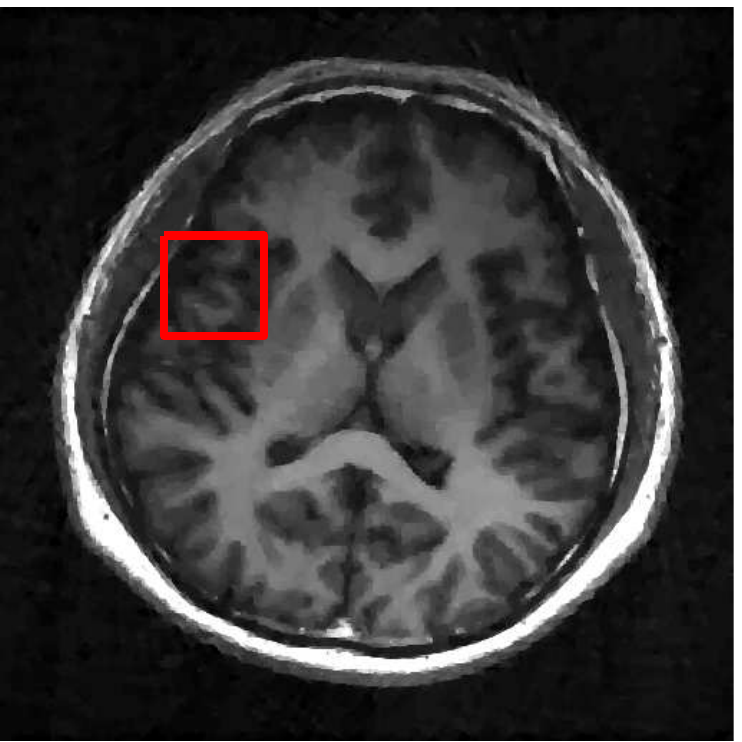}&
\includegraphics[width=.23\textwidth]{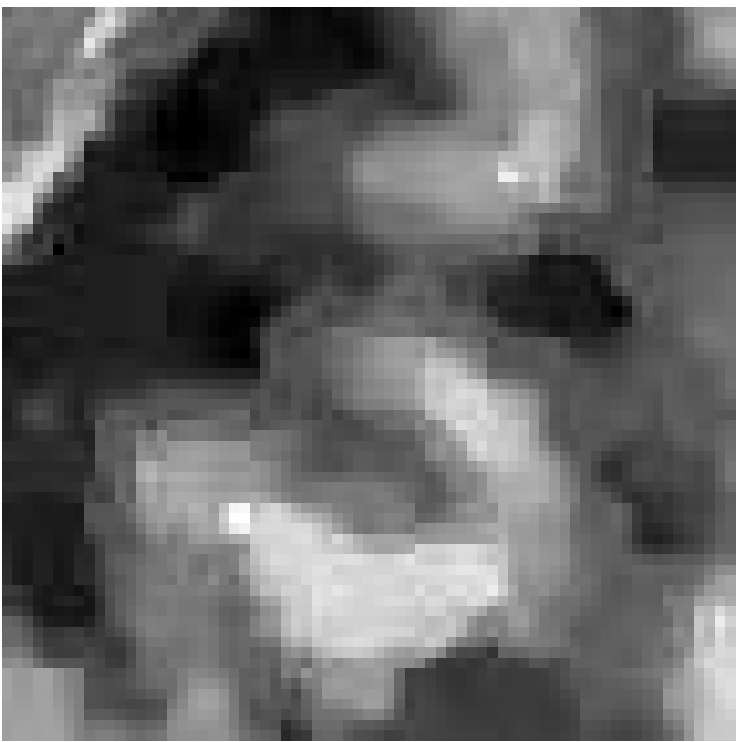}
\end{tabular}
\end{center}
\caption{Reconstructed brain MR image. First row from left to right: ground truth, close-up of ground truth, our result, close-up of the result.
Second row from left to right: result obtained by RecPF, close-up of RecPF result, result by EdgeCS, close-up of EdgeCS result.}
\label{fig:test1}
\end{figure*}

\begin{figure*}
\begin{center}
\includegraphics[width=.32\textwidth]{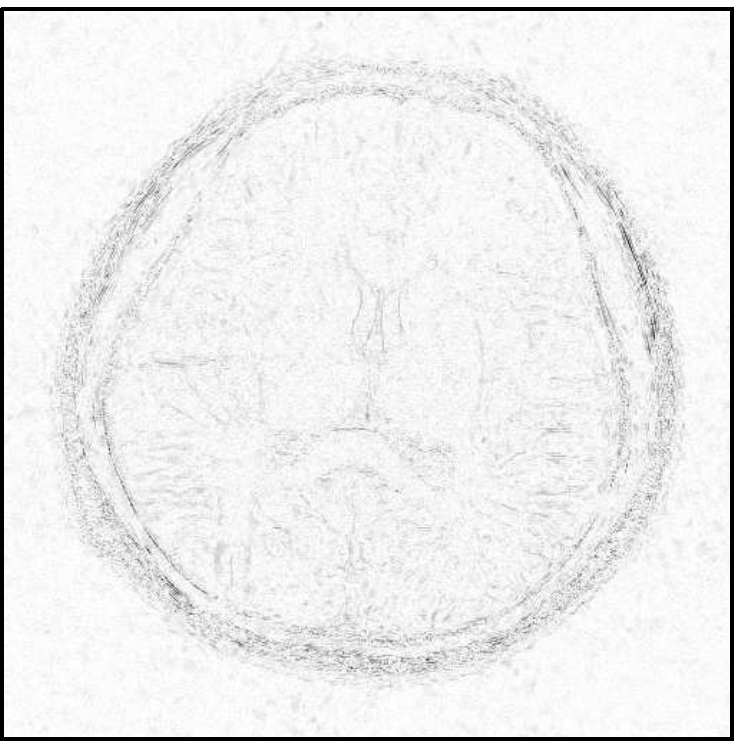}
\includegraphics[width=.32\textwidth]{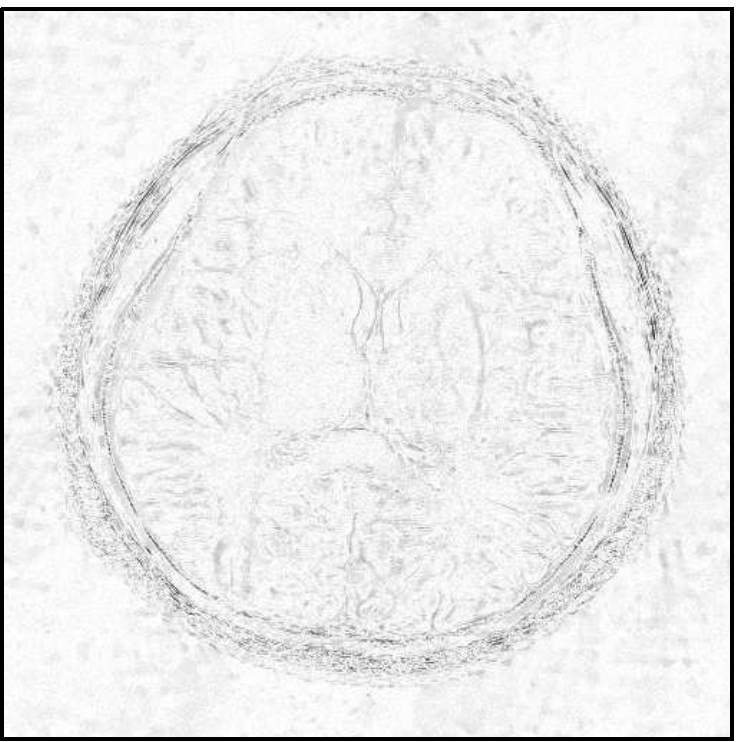}
\includegraphics[width=.32\textwidth]{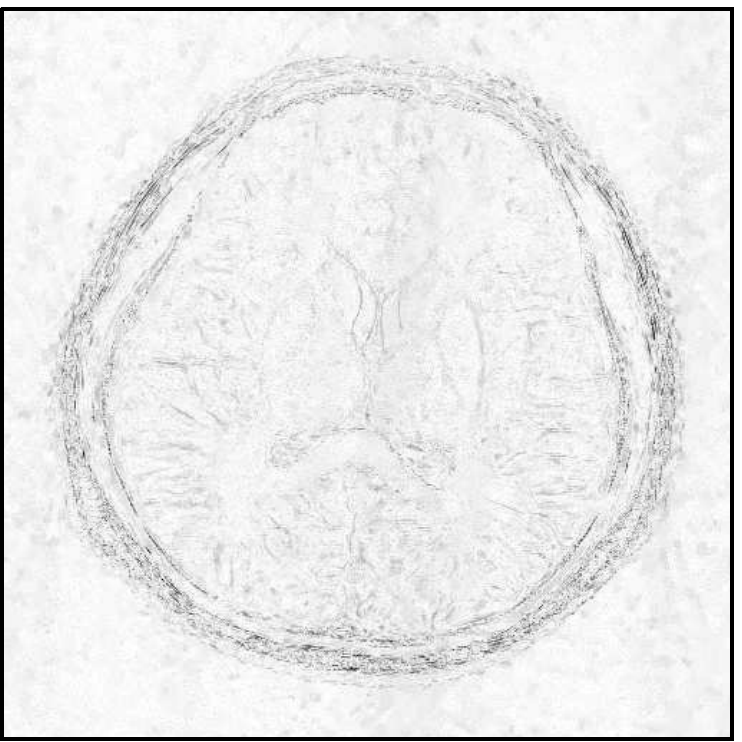}
\end{center}
\caption{Residual maps of the three method results with the ground truth. For better visualization, we inverted the grayscale. From left to right: proposed, RecPF, EdgeCS, the relative errors are listed respectively as: 12.00\%, 14.68\%, 15.81\%.}
\label{fig:test1dif}
\end{figure*}

In the first example, we look at the simulated noise-free spectral measurements of a $512\times 512$ brain MR image downloaded from BrainWeb \url{https://brainweb.bic.mni.mcgill.ca/brainweb/}. The ground truth image has inhomogenous contrasts in different areas, especially in the gray matter and cerebrospinal fluid. We tested the proposed GeoCS algorithm, RecPF and EdgeCS with 40 radial sampling lines, namely 8.79\% sampling rate. We show the results in Fig. \ref{fig:test1} and zoom in one small patch for better visual comparison. It's apparent that the image produced by GeoCS has better quality than the others. RecPF sort of oversmooths the whole image, and EdgeCS is able to detect the edges while losing some gradual transition between smooth areas and boundaries. To further compare three results, we take the difference between the ground truth and the reconstructed image for each method and display the inverted residue images in Fig. \ref{fig:test1dif}. It's clear that our proposed algorithm suppresses the error more evenly inside the skull. The three approaches are also compared as the sampling rate changes. The quantitative comparison listed in Table \ref{table:test1err} shows that the proposed GeoCS consistently outperforms the other methods.

\begin{table}
\centering
\begin{tabular}{l|c|c|c|c|c|c|c|c}
\hline\hline
\multirow{2}{*}{Sampling rate} & \multicolumn{2}{|c|}{8.79\%} &\multicolumn{2}{|c|}{10.64\%} & \multicolumn{2}{|c|}{12.92\%} & \multicolumn{2}{|c}{14.74\%} \\ \cline{2-9}
& RelErr & SNR & RelErr & SNR & RelErr & SNR & RelErr & SNR \\ \hline
Proposed &\textbf{ 0.1200} & \textbf{15.90} & \textbf{0.1016} & \textbf{17.34} & \textbf{0.0874} &\textbf{18.65} & \textbf{0.0797} &\textbf{19.46}\\
RecPF & 0.1468 & 14.16 & 0.1273 & 15.40 & 0.1097 & 16.68 & 0.1011 & 17.39 \\
EdgeCS & 0.1581 & 14.06 & 0.1378 & 15.34 & 0.1208 & 16.56 & 0.1111 & 17.35\\ \hline\hline
\end{tabular}
\caption{Relative error and SNR comparisons for brain MRI reconstruction}
\label{table:test1err}
\end{table}


\subsection{Example 2}
\begin{figure*}
\begin{center}
\begin{tabular}{cccc}
\includegraphics[width=.23\textwidth]{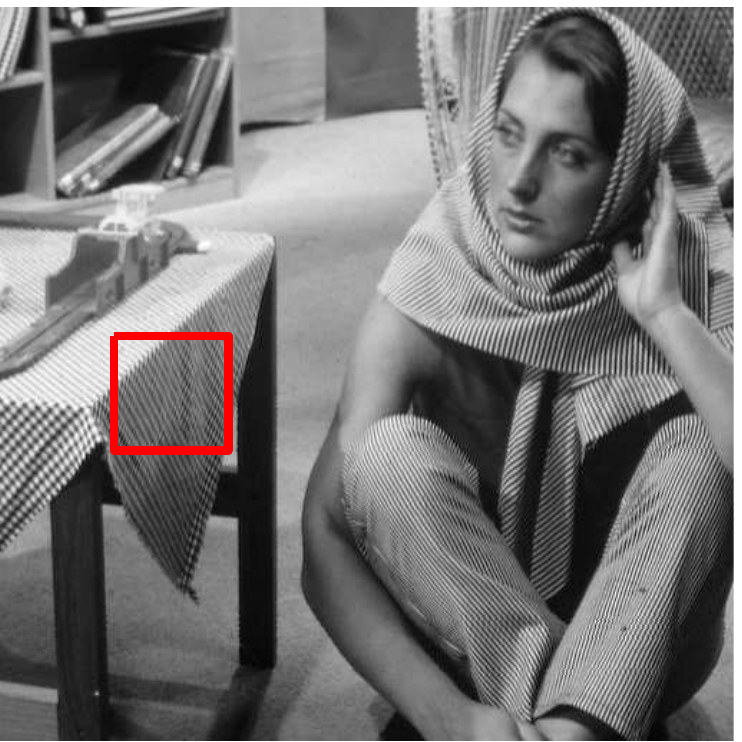}&
\includegraphics[width=.23\textwidth]{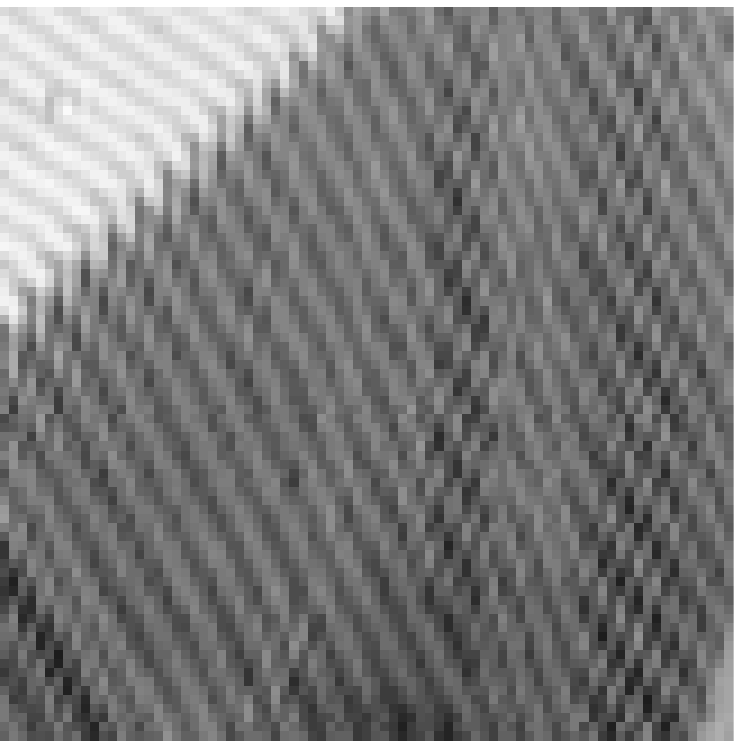}&
\includegraphics[width=.23\textwidth]{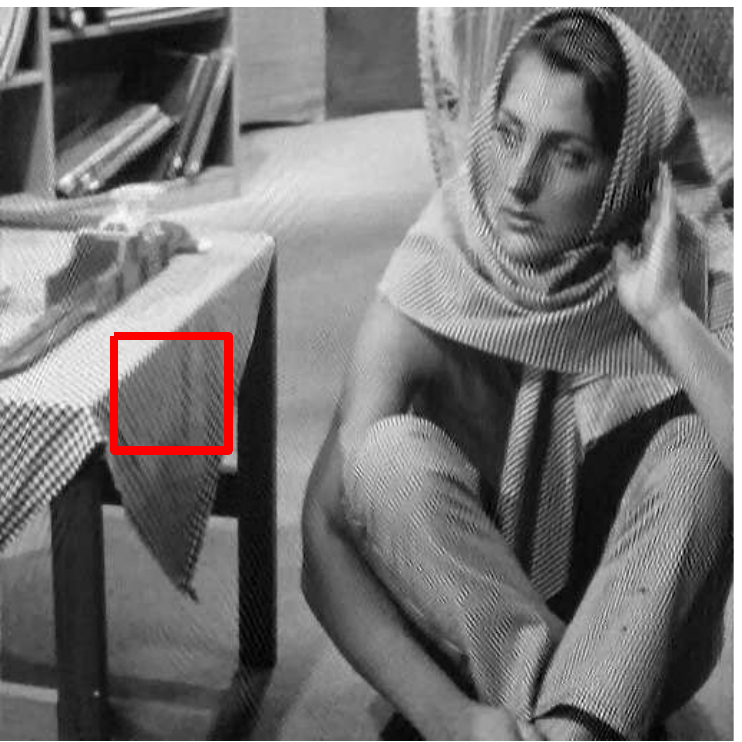}&
\includegraphics[width=.23\textwidth]{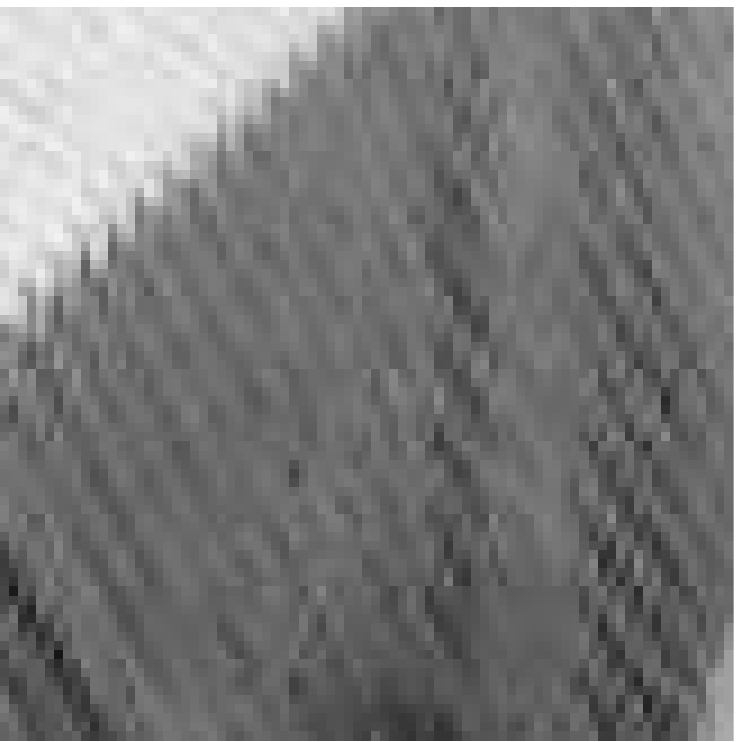}\\
\includegraphics[width=.23\textwidth]{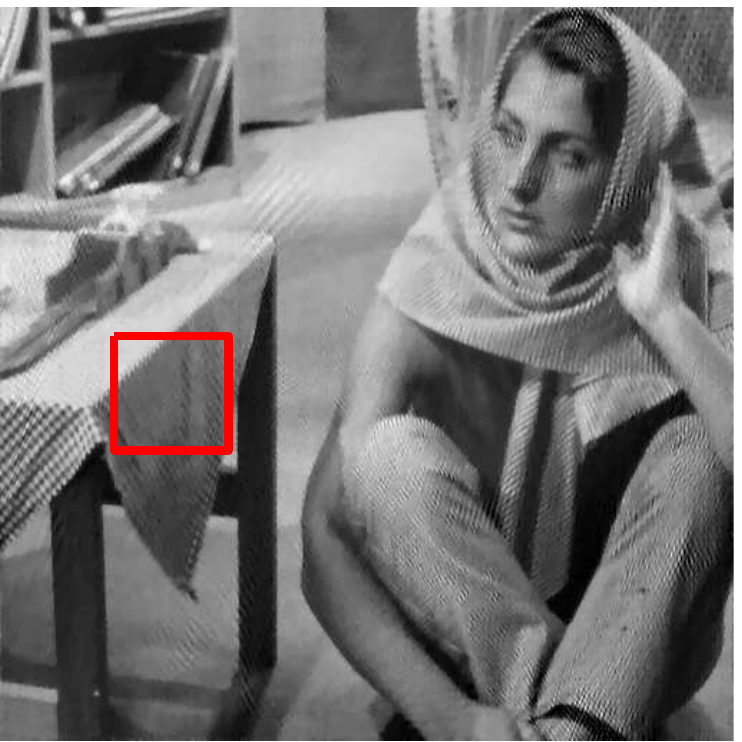}&
\includegraphics[width=.23\textwidth]{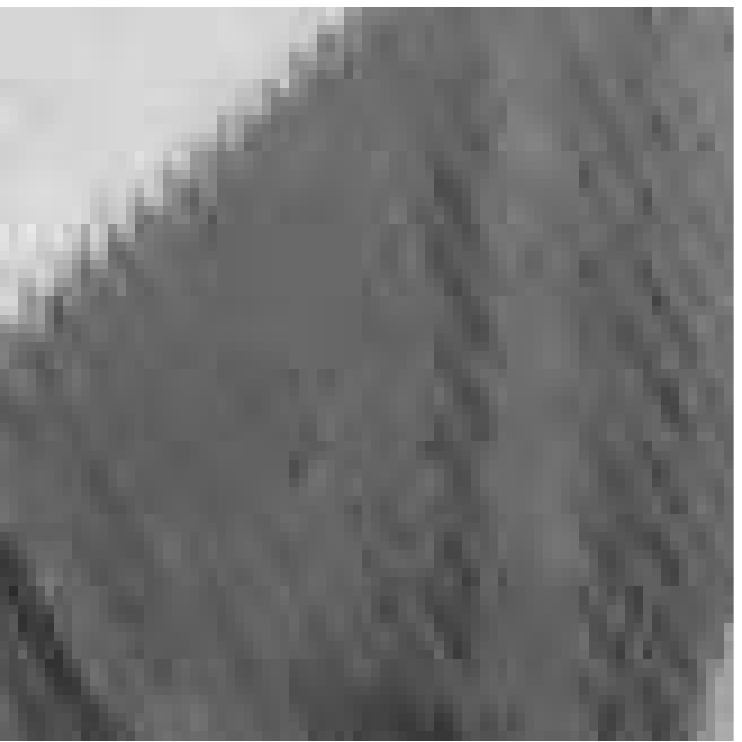}&
\includegraphics[width=.23\textwidth]{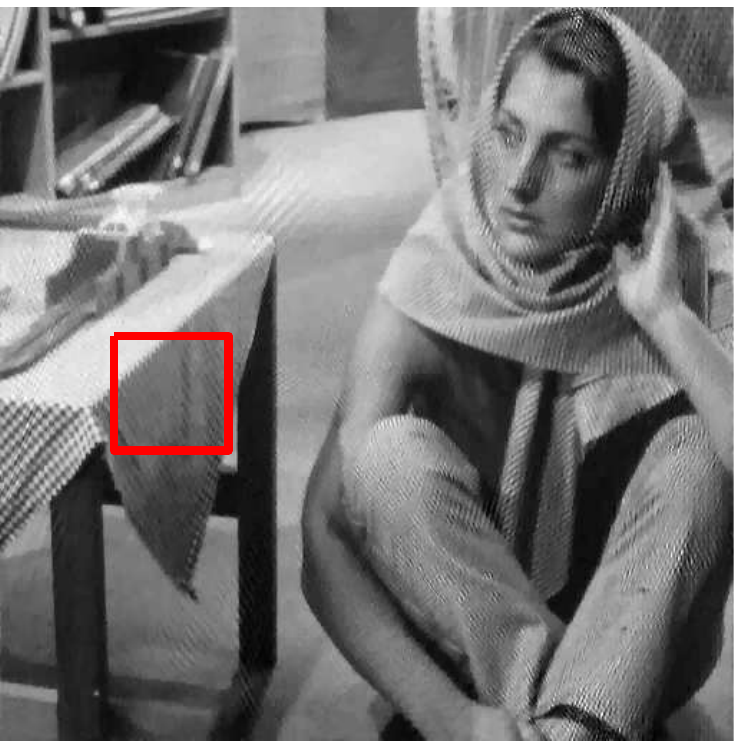}&
\includegraphics[width=.23\textwidth]{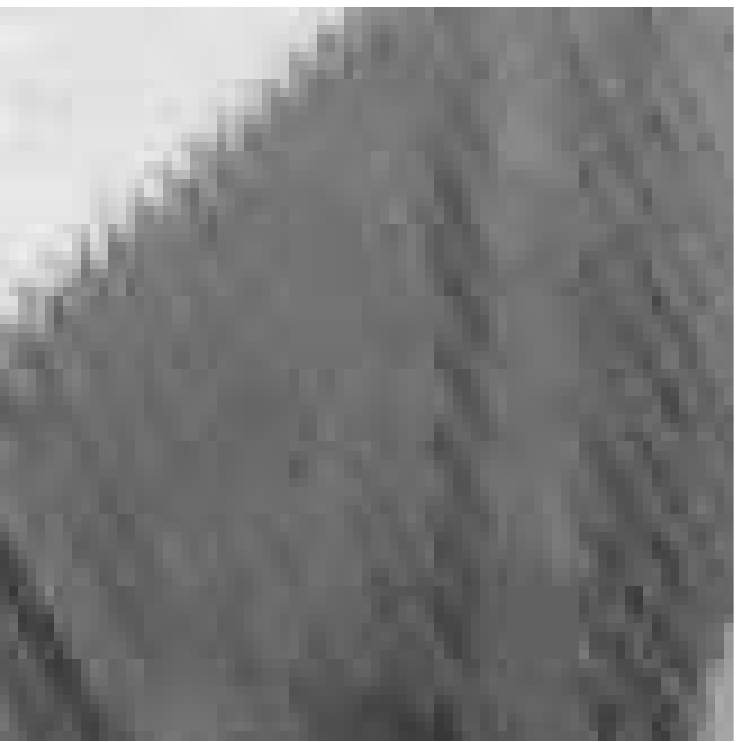}
\end{tabular}
\end{center}
\caption{Reconstructed barbara image. First row from left to right: accurate barbara image, cropped barbara image, our result, close-up of our result.
Second row from left to right: result by RecPF, close-up of RecPF result, result by EdgeCS, close-up of EdgeCS result.}
\label{fig:test2}
\end{figure*}

\begin{figure*}
\begin{center}
\includegraphics[width=.32\textwidth]{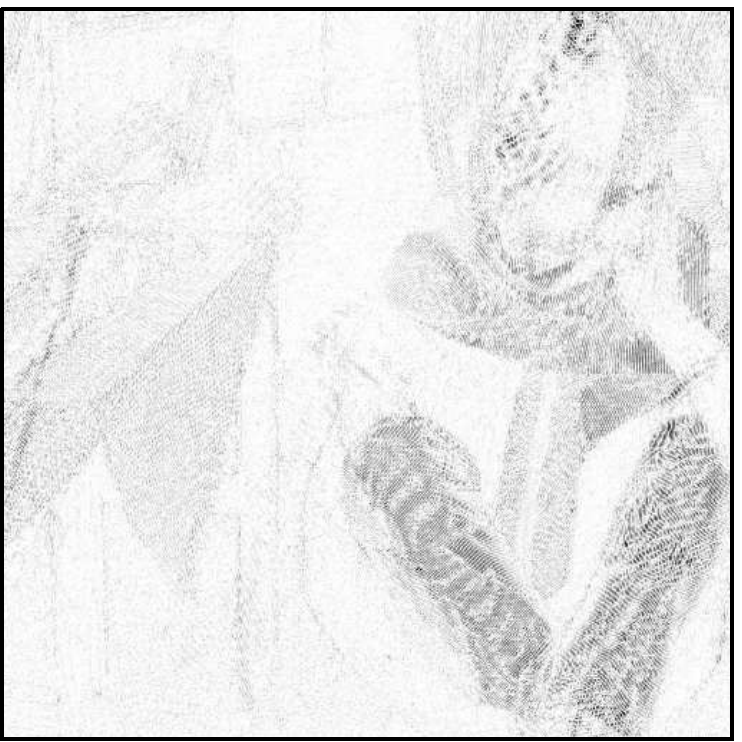}
\includegraphics[width=.32\textwidth]{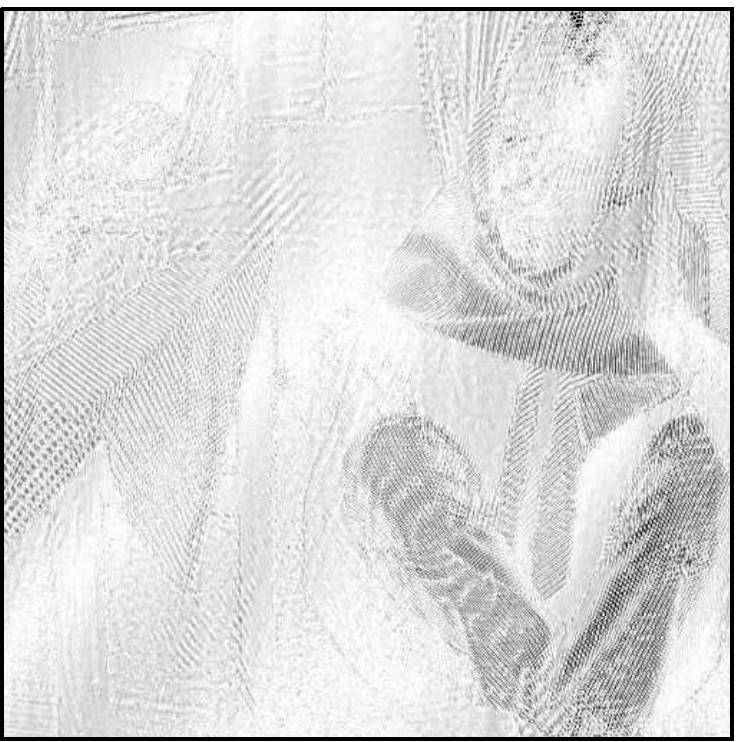}
\includegraphics[width=.32\textwidth]{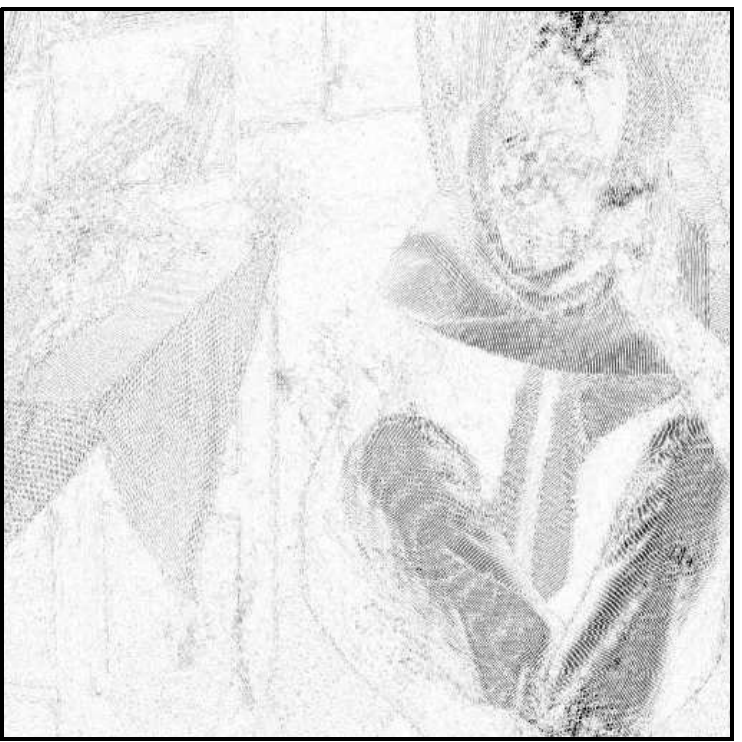}
\end{center}
\caption{From left to right: proposed, RecPF, EdgeCS. Relative error: 9.64\%, 14.01\%, 12.53\%.}
\label{fig:test2dif}
\end{figure*}

In this example, we show the benefits of GeoCS on reconstructing texture images. The test image is $512\times512$ Barbara image which has various texture patterns and a lot of details. This image requires relatively higher sampling rate to get an ideal recovery and standard edge detection algorithm may even fail to get accurate edges. The results obtained by GeoCS, RecPF and EdgeCS with 100 radial sampling lines (sampling rate 20.87\%) are listed in Fig. \ref{fig:test2}, where we zoomed in one patch of table cloth. Our proposed method is able to recover largely the directional textures while the other two methods get blurry textures. The inverted residue images are listed in Fig. \ref{fig:test2dif}. The consistent performance is illustrated in Table \ref{table:test2err} using different sampling rates. In this example, shearlet transform plays an important role in preserving the structures in different directions and thereby the textures. Weighted total variation further corrects the smooth areas which were over-texturized.

\begin{table}
\centering
\begin{tabular}{l|c|c|c|c|c|c|c|c}
\hline\hline
\multirow{2}{*}{Sampling rate} & \multicolumn{2}{|c|}{8.79\%} & \multicolumn{2}{|c|}{12.92\%} & \multicolumn{2}{|c|}{16.94\%} & \multicolumn{2}{|c}{20.87\%}\\ \cline{2-9}
& RelErr & SNR & RelErr & SNR & RelErr & SNR & RelErr & SNR \\ \hline
Proposed &\textbf{0.1429} & \textbf{10.16} & \textbf{0.1248} & \textbf{11.33} & \textbf{0.1104}& \textbf{12.39} & \textbf{0.0964} &\textbf{13.58} \\
RecPF  & 0.1700 & 8.66 & 0.1552 & 9.45  & 0.1471 & 9.92  & 0.1401 & 10.36\\
EdgeCS & 0.1574 & 9.69 & 0.1412 & 10.72 & 0.1323 & 11.35 & 0.1253 & 11.92\\ \hline\hline
\end{tabular}
\caption{Relative error and SNR comparisons for barbara image reconstruction}
\label{table:test2err}
\end{table}


\subsection{Example 3}
\begin{figure*}
\begin{center}
\begin{tabular}{cccc}
\includegraphics[width=.23\textwidth]{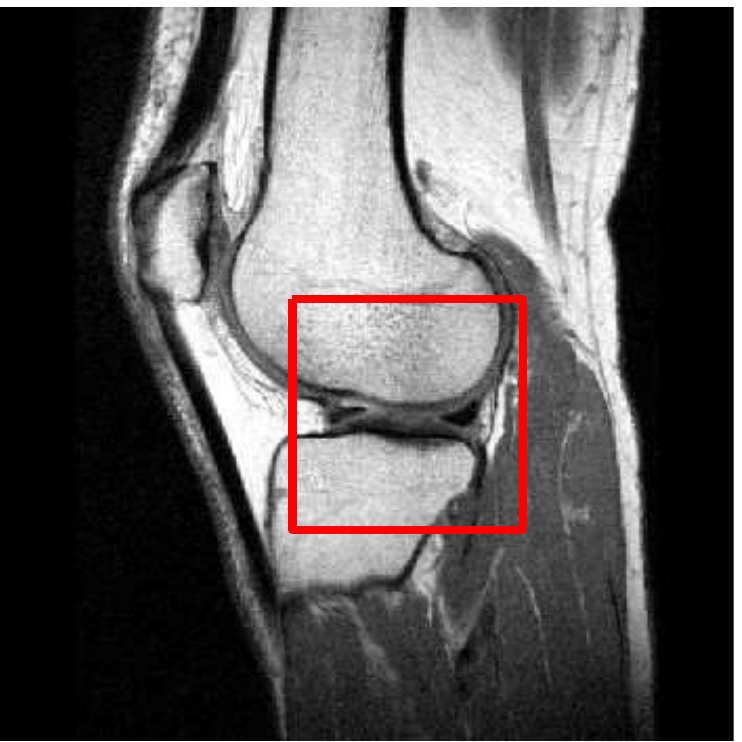}&
\includegraphics[width=.23\textwidth]{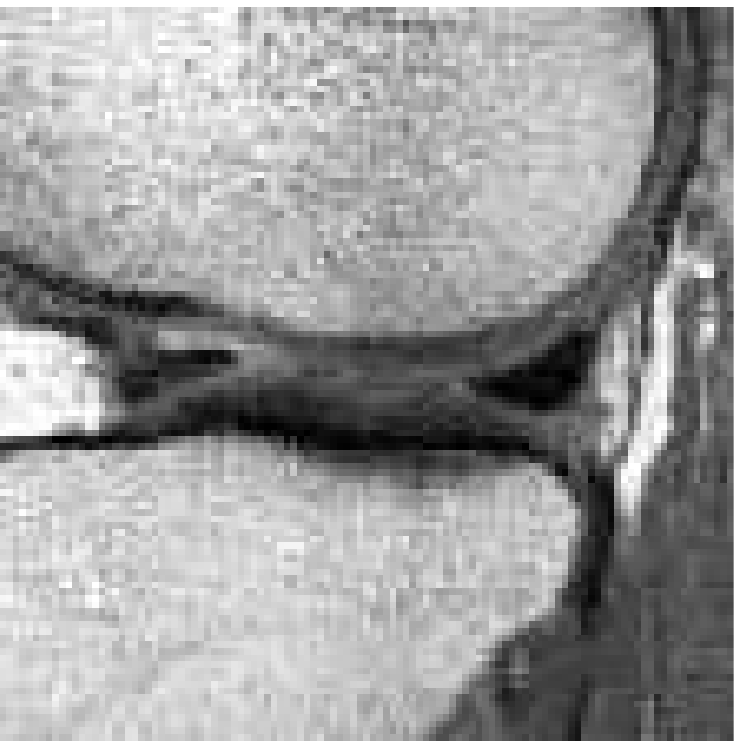}&
\includegraphics[width=.23\textwidth]{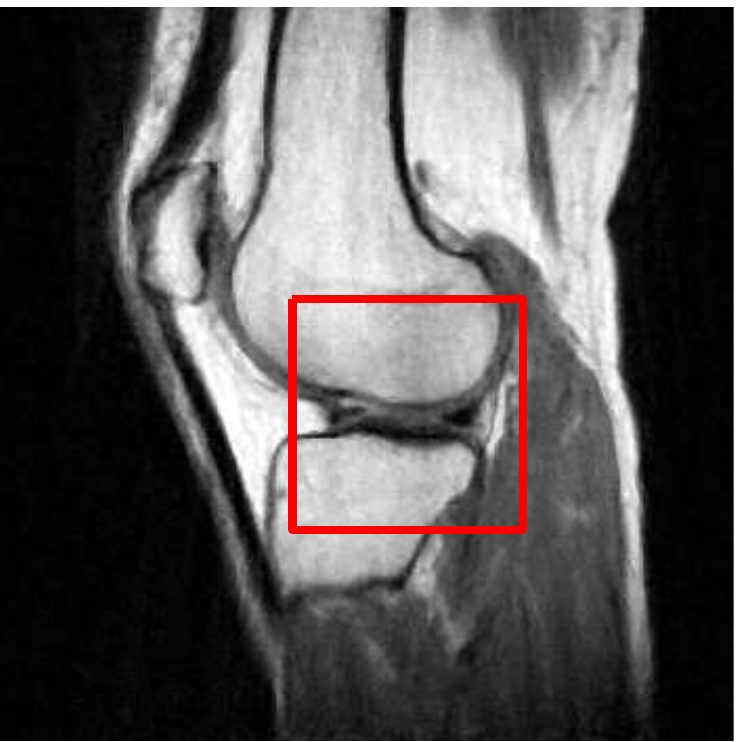}&
\includegraphics[width=.23\textwidth]{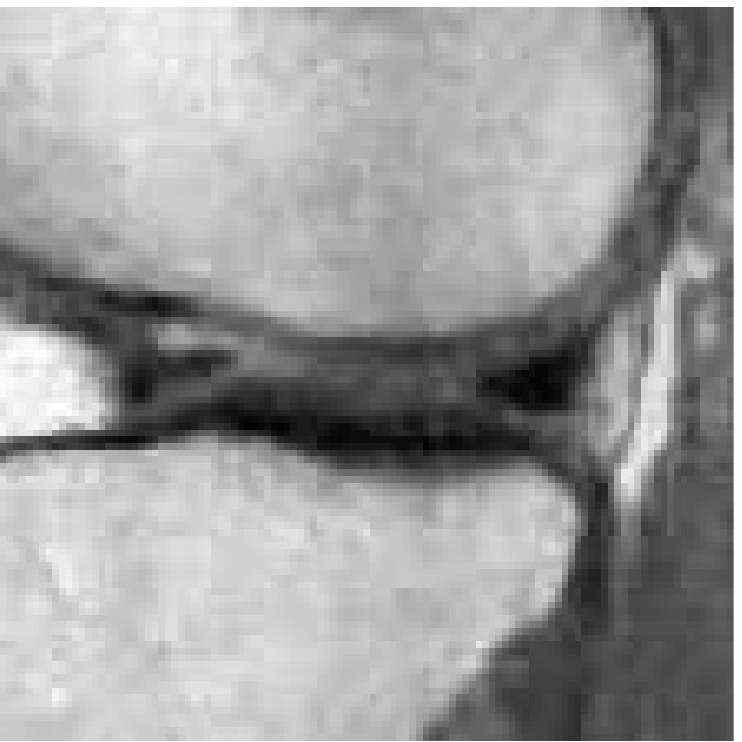}\\
\includegraphics[width=.23\textwidth]{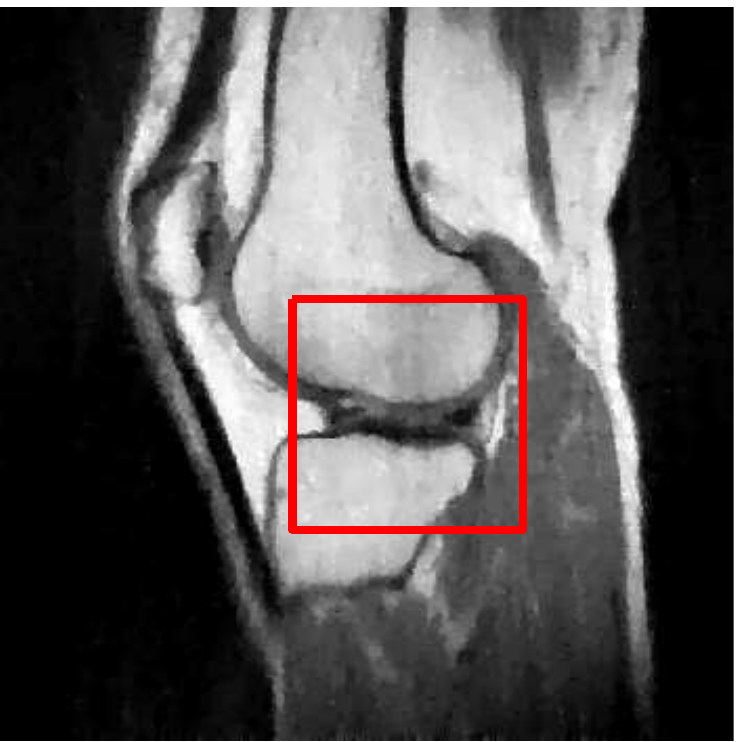}&
\includegraphics[width=.23\textwidth]{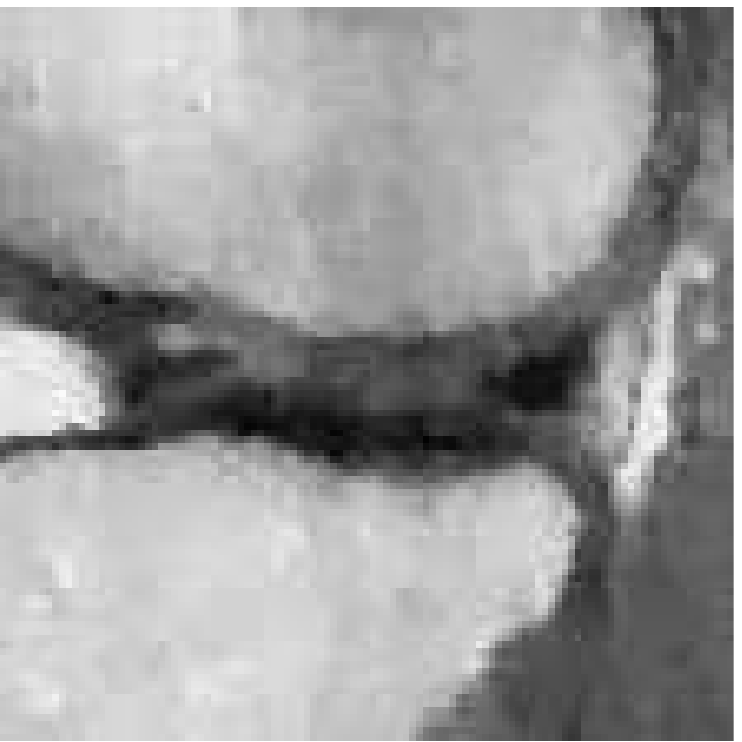}&
\includegraphics[width=.23\textwidth]{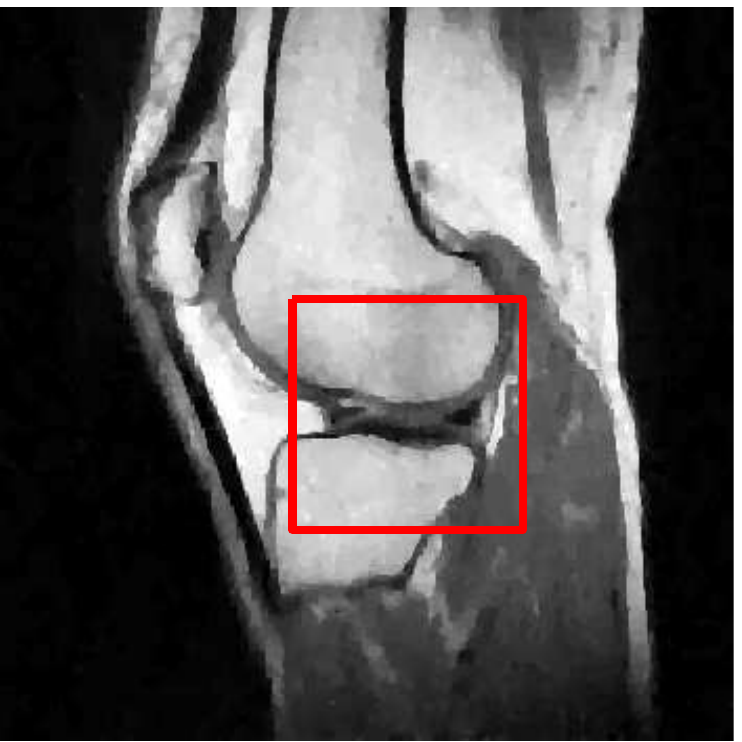}&
\includegraphics[width=.23\textwidth]{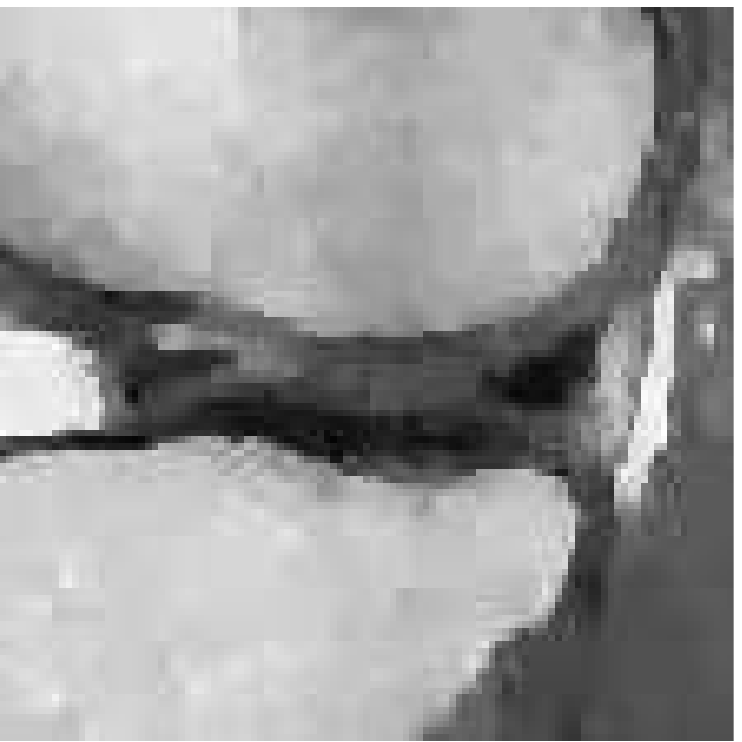}
\end{tabular}
\end{center}
\caption{Recovered knee image from noisy data. Top row from left to right: ground truth, our result. Bottom row from left to right: result by RecPF and EdgeCS.}
\label{fig:test3}
\end{figure*}

\begin{figure*}
\begin{center}
\includegraphics[width=.32\textwidth]{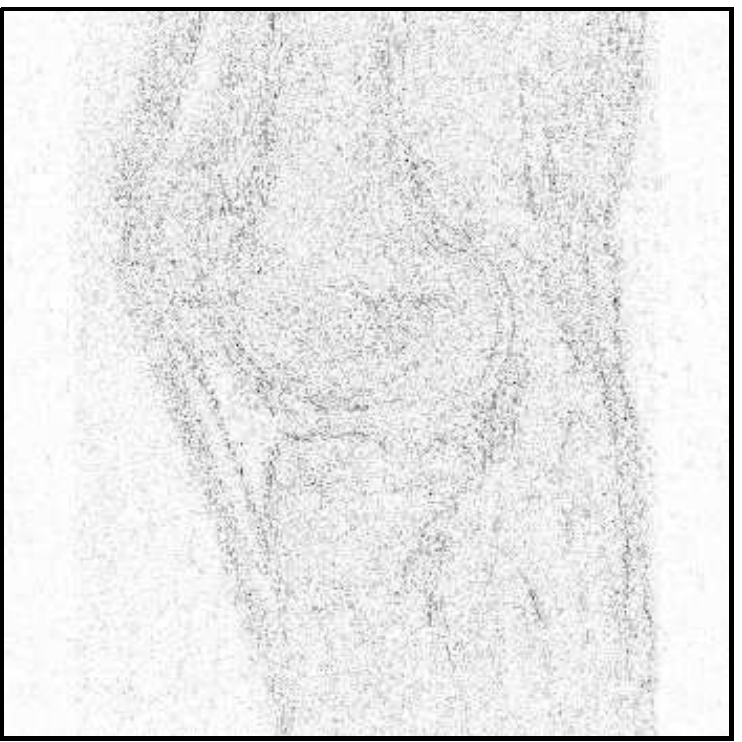}
\includegraphics[width=.32\textwidth]{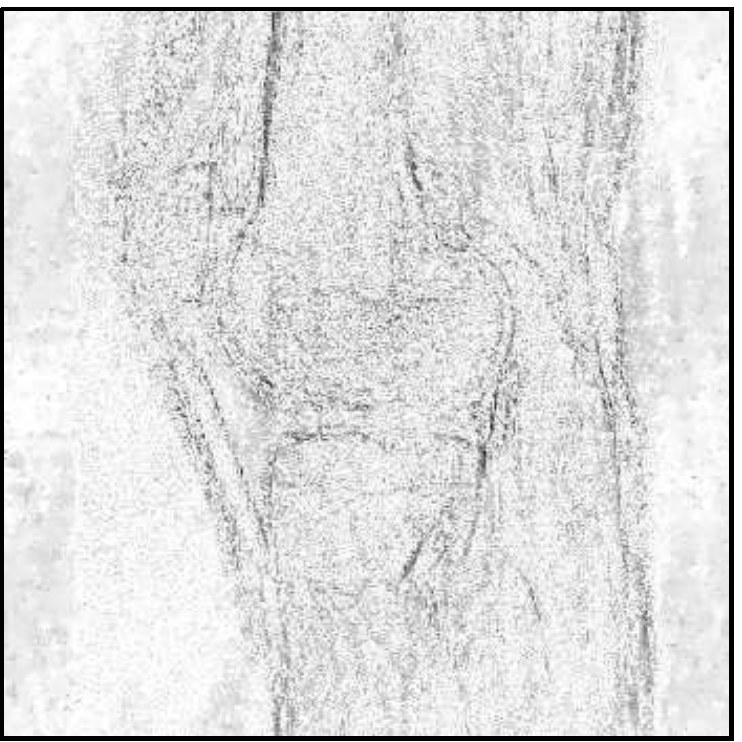}
\includegraphics[width=.32\textwidth]{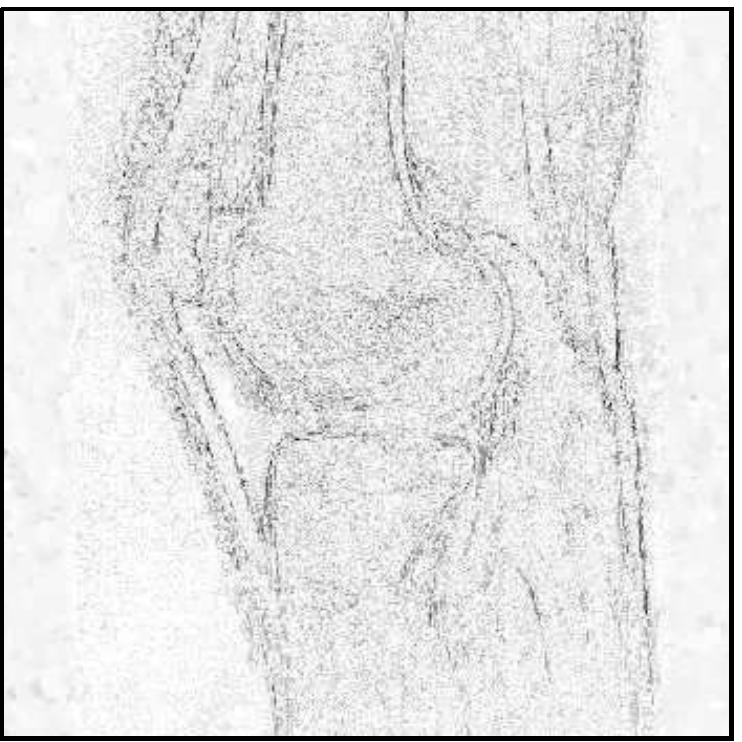}
\end{center}
\caption{From left to right: proposed, RecPF, EdgeCS.
Relative error: 9.52\%, 11.71\%, 10.35\%.}
\label{fig:test3dif}
\end{figure*}

Our last test image is a T1 weighted MR image of the knee showing femur, patella, tibia and menisci from \url{http://www.mr-tip.com/}. We first added zero-mean complex Gaussian noise $\sigma = 10$ to the spectral data sampled by 40 radial lines (sampling rate 12.71\%). The recovered images and their associated enlarged patches given by GeoCS, RecPF and EdgeCS are shown respectively in Fig. \ref{fig:test3}. Visually our result is more natural in the bones and junctions and is closer to the ground truth than those given by the other two methods. It is worth noting that staircase artifacts are significantly reduced in our result while they are quite obvious in the RecPF and EdgeCS results especially in the smooth bone regions. For better visual comparison, we show the residue images in Fig. \ref{fig:test3dif}.

By fixing the sampling rate as 12.71\%, we also perform the comparisons on spectral data with different noise levels $\sigma=5,10,15,20$. In the case of large noise level, it is better to adjust the regularization parameters accordingly. It is also true in RecPF and EdgeCS algorithms. From Table \ref{table:test3err}, we observe that our proposed algorithm is robust to the noise and produces more accurate reconstructed images than the other methods. All the results are obtained under the optimal parameter settings for each method.

\begin{table}
\centering
\begin{tabular}{l|c|c|c|c|c|c|c|c}
\hline\hline
\multirow{2}{*}{\phantom{xx}$\sigma$} & \multicolumn{2}{|c|}{5} & \multicolumn{2}{|c|}{10} & \multicolumn{2}{|c|}{15} & \multicolumn{2}{|c}{20}\\ \cline{2-9}
& RelErr & SNR & RelErr & SNR & RelErr & SNR & RelErr & SNR \\ \hline
Proposed     &\textbf{0.0842}& \textbf{18.50} &\textbf{0.0952}& \textbf{17.44} & \textbf{0.1006}& \textbf{16.95} & \textbf{0.1095} & \textbf{16.22}\\
RecPF            &0.1064& 16.47 &0.1171& 15.64 &0.1283 & 14.85 &0.1389 & 14.16 \\
EdgeCS           &0.0992& 17.08 &0.1035& 16.71 &0.1118 & 16.04 &0.1194 & 15.47 \\ \hline\hline
\end{tabular}
\caption{Relative error and SNR comparisons for the noisy knee MRI reconstruction}
\label{table:test3err}
\end{table}


\section{Conclusion and remarks}\label{conclusion}
We proposed a two-stage compressive sensing image reconstruction algorithm based on shearlet transform and weighted TV. The first stage is to use standard $\TV$-$L_1$-$L_2$ model with shearlet transform to get an initial guess for the underlying image of interest. Geometric information extracted from this guess serves as an initial a priori in weighted $\TV$-$L_1$-$L_2$ model to further enhance the reconstruction accuracy. This kind of geometric information extraction and image reconstruction are alternated in a mutually beneficial fashion until it converges. Replacing the conventional wavelet transform with shearlet transform, the model is able to promote the signal sparsity, and preserve multiple directional features better during the recovery. The spatially variant weights associated to TV plays an important part in preserving sharp edges while reducing staircase effects of TV. The minimization problem is solved by split Bregman which divides one complicated optimization problem with nondifferentiable terms into three sets of subproblem, each of which has closed-form solutions. Convergence of the algorithm is guaranteed under mild conditions. The proposed approach is compared with two recent related work. Numerical experiments show the consistent overwhelming advantages of our algorithm.

The proposed approach GeoCS is better than RecPF and EdgeCS in reconstructing complicated piecewise smooth images. However, as for piecewise constant images, it is sufficient to apply one-stage methods. Moreover, by adapting the weights to the spatially variant gradients along with two-stage reweighting scheme, GeoCS integrates more reliable geometric prior to the reconstruction than EdgeCS. Our extensive experience shows that the more accurate geometric information is obtained during the algorithm, the better the overall scheme will perform. A high-quality result of Stage I will speed up the convergence of Stage II and potentially ease the parameter tuning.
Nevertheless, it is possible that the extracted geometric information is not reliable at all in case of extremely insufficient samples or excessive noise and thereby GeoCS might fail. There is still room to study how to extract much more reliable geometric information from noisy incomplete measurements and how to efficiently utilize them. Furthermore, some recent acceleration techniques, e.g., Nesterov's accelerated gradient descent\cite{Nesterov:83}, can be applied to speed up the convergence.

\begin{acknowledgement}
The authors would like to thank the Research Collaboration Workshop for Women in Data Science and Mathematics held at ICERM during July 29-August 2, 2019. Qin is supported by the NSF grant DMS-1941197, and Guo is supported by the NSF grant DMS-1521582.
\end{acknowledgement}

\bibliographystyle{unsrt}
\bibliography{ref}

\end{document}